\documentclass[a4paper,11pt,twoside]{article}
\usepackage{amssymb,amsmath,latexsym,geometry,fancyhdr,lineno,hyperref,titletoc}
\contentsmargin{0pt}

\dottedcontents{section}[2em]{\vspace{-1mm}\small}{2em}{0pt}
\dottedcontents{subsection}[5em]{\vspace{-1mm}\small}{3em}{5pt}

\hypersetup{colorlinks,linkcolor= blue,citecolor=blue}
\newtheorem{theorem}{Theorem}[section]
\newtheorem{lemma}[theorem]{Lemma}

\newtheorem{remark}{Remark}[section]
\newtheorem{definition}{Definition}[section]
\newtheorem{corollary}[theorem]{Corollary}
\def\proof{\mbox {\textbf{Proof.}~~}}
\numberwithin{equation}{section}
\allowdisplaybreaks[4]
\begin{document}
\title{{\bf\Large Nonlinear scalar field equations with general nonlinearity}}
\author{\\
{ \textbf{\normalsize Louis Jeanjean}}\footnote{For L. Jeanjean this work has been carried out in the framework of the project NONLOCAL (ANR-14-CE25-0013) funded by the French National Research Agency (ANR).}\\
{\it\small Laboratoire de Math\'{e}matiques (CNRS UMR 6623)}\\
{\it\small Universit\'{e} de Bourgogne Franche-Comt\'{e}}\\
{\it\small Besan\c{c}on, 25030, France}\\
{\it\small e-mail: louis.jeanjean@univ-fcomte.fr}\\
\\
{ \textbf{\normalsize Sheng-Sen Lu}}\footnote{S.-S. Lu acknowledges the support of the NSF of China (NSFC-11771324), of the China Scholarship Council (CSC-201706250149) and the hospitality the Laboratoire de Math\'{e}matiques (CNRS UMR 6623), Universit\'{e} de Bourgogne Franche-Comt\'{e}.}\\
{\it\small Center for Applied Mathematics, Tianjin University}\\
{\it\small Tianjin, 300072, PR China}\\
{\it\small e-mail: sslu@tju.edu.cn}}
\date{}
\maketitle
{\bf\normalsize Abstract.} {\small Consider the nonlinear scalar field equation
\begin{equation}\label{a1}
  -\Delta{u}= f(u)\quad\text{in}~\mathbb{R}^N,\qquad u\in H^1(\mathbb{R}^N),
\end{equation}
where $N\geq3$ and $f$ satisfies the general Berestycki-Lions conditions. We are interested in the existence of positive ground states, of nonradial solutions and in the multiplicity of radial and nonradial solutions. Very recently Mederski \cite{Me17} made a major advance in that direction through the development, in an abstract setting, of a new critical point theory for constrained functionals.  In this paper we propose an alternative, more elementary approach, which permits to recover Mederski's results on \eqref{a1}. The keys to our approach are an extension to the symmetric mountain pass setting of the monotonicity trick, and a new decomposition result for bounded Palais-Smale sequences.}




\pagestyle{fancy}
\fancyhead{} 
\fancyfoot{} 
\renewcommand{\headrulewidth}{0pt}
\renewcommand{\footrulewidth}{0pt}
\fancyhead[CE]{ Nonlinear scalar field equations with general nonlinearity}
\fancyhead[CO]{ Nonlinear scalar field equations with general nonlinearity}
\fancyfoot[C]{\thepage}
\tableofcontents

\section{Introduction}\label{sect:introduction}
We consider nonlinear scalar field equations
\begin{linenomath*}
  \begin{equation}\label{mainproblem}
    -\Delta{u}= f(u)\quad\text{in}~\mathbb{R}^N,\qquad u\in H^1(\mathbb{R}^N),
  \end{equation}
\end{linenomath*}
where $N\geq3$ and $f$ satisfies the assumptions stated below:
\begin{itemize}
  \item[$(f1)$] $f:\mathbb{R}\to\mathbb{R}$ is continuous and odd.
  \item[$(f2)$] $-\infty<\liminf_{t\rightarrow0}f(t)/t\leq\limsup_{t\rightarrow0}f(t)/t<0$.
  \item[$(f3)$] $\lim_{t\rightarrow\infty}f(t)/|t|^{\frac{N+2}{N-2}}=0$.
  \item[$(f4)$] There exists $\zeta>0$ such that $F(\zeta)>0$, where $F(t):=\int^t_0f(s)ds$ for $t\in\mathbb{R}$.
\end{itemize}

In the fundamental papers \cite{Be83-1,Be83-2}, Berestycki and Lions introduced the assumptions $(f1)-(f4)$ for the first time when dealing with \eqref{mainproblem}. The main feature of these assumptions is that they are almost necessary to get a nontrivial solution to \eqref{mainproblem}, see however Remark \ref{extensionf3} below.  It turns out that these assumptions are also sufficient. Indeed, with the aid of variational methods, by studying certain constrained problems, Berestycki-Lions showed the existence of a \emph{ground state solution} which is positive and radially symmetric in \cite{Be83-1}, and obtained infinitely many radial solutions in \cite{Be83-2}.

After these two papers by Berestycki and Lions, still under the assumptions $(f1)-(f4)$, several advances in the understanding of the set of solutions to \eqref{mainproblem} (including but not limited to those listed below) were made in subsequent works. From now and throughout the paper, Problem \eqref{mainproblem} will refer to equation \eqref{mainproblem} considered under the assumptions $(f1)-(f4)$.

An important observation on Problem \eqref{mainproblem} is pointed out in \cite{Je03}. In that paper, the authors showed that the associated energy functional
\begin{linenomath*}
\begin{equation*}
J(u):=\frac{1}{2}\int_{\mathbb{R}^N}|\nabla u|^2dx-\int_{\mathbb{R}^N}F(u)dx
\end{equation*}
\end{linenomath*}
has a mountain pass geometry, that is,
\begin{linenomath*}
\begin{equation*}
  c_{mp}:=\underset{\gamma\in\Gamma}{\inf}\underset{t\in[0,1]}{\max}J(\gamma(t))>0,
\end{equation*}
\end{linenomath*}
where $\Gamma:=\left\{\gamma\in C([0,1], H^1(\mathbb{R}^N))~|~\gamma(0)=0,J(\gamma(1))<0\right\}$. It was also proved that the ground state solutions are actually mountain pass solutions, which are considered as structurally stable. This fact turns out to be very useful in the studies of the corresponding singular perturbation problems and non-autonomous cases, especially when one tries to relax the more restricted conditions assumed on the nonlinearity to the almost optimal ones like $(f2)-(f4)$, see, e.g., \cite{By10,By07-1,By12,Je05}.

We also would like to mention the work \cite{By09}, in which it was shown that any ground state solution of Problem \eqref{mainproblem} is radially symmetric (up to a translation), has a constant sign and is monotone with respect to the radial variable.

In a more recent paper \cite{Hi10}, Hirata, Ikoma and Tanaka revisited Problem \eqref{mainproblem} in $H^1_{\mathcal{O}}(\mathbb{R}^N)$ the subspace of radially symmetric functions of $H^1(\mathbb{R}^N)$, and showed further that the functional $J$ has a symmetric mountain pass geometry. By using the mountain pass and symmetric mountain pass approaches developed in \cite{Am73}, they managed to find a positive \emph{radial ground state solution} (namely, a nontrivial radial solution minimizing $J$ among all the nontrivial radial solutions) and infinitely many radial solutions through the unconstrained functional $J$. This is in contrast to the situation in \cite{Be83-1,Be83-2} where the solutions were constructed through certain constrained problems. The core of the proof developed in \cite{Hi10} is the use of a suitable extended functional on the augmented space $\mathbb{R}\times H^1_\mathcal{O}(\mathbb{R}^N)$. This technique of \emph{adding one dimension of space} was first introduced in \cite{Je97} to deal with a nonlinear eigenvalue problem, and we refer readers to \cite{Az11,Ba13,Ch14,Cu15,Hi18,Lu16,Mo15} for its recent applications to various problems.

The most recent advance on \eqref{mainproblem}, assuming just $(f1)-(f4)$, was made by Mederski in \cite{Me17} where the existence and multiplicity of \emph{nonradial solutions} to \eqref{mainproblem}, were established for the first time. More precisely, Mederski found at least one nonradial solution for any $N\geq4$, and showed the existence of infinitely many distinct nonradial solutions if in addition $N\neq5$. These results give a partial positive answer to a question which was posed by Berestycki and Lions (see \cite[Section 10.8]{Be83-2}) and had been open for more than thirty years. The proofs in \cite{Me17} are based on a new constrained approach, developed in an abstract setting, applied to treat \eqref{mainproblem}, see \cite{Me17} for more details.   As another application of his approach, Mederski gave a new proof of the existence of a ground state  for Problem \eqref{mainproblem}.

Let us now present the main results of this paper. As it will be clear these results are not new but to derive them we propose a new approach which we believe has its own interest. Theorems \ref{theorem:groundstate} and \ref{theorem:radialsolutions} stated below which concern ground states and radial solutions respectively, are well known since the papers \cite{Be83-1,Be83-2}.
\begin{theorem}\label{theorem:groundstate}
 Assume that $(f1)-(f4)$ hold, then \eqref{mainproblem} has a positive ground state solution.
\end{theorem}
\begin{remark}
As one will see, Theorem \ref{theorem:groundstate} will be proved by a mountain pass argument in $H^1(\mathbb{R}^N)$. Even though
it was already shown in \cite{Je03} that the functional $J$ has a mountain pass geometry and that the mountain pass value $c_{mp}$ corresponds to the infimum of the nonzero critical levels of $J$, a direct proof as we present here seemed to be missing.
\end{remark}
\begin{theorem}\label{theorem:radialsolutions}
 Assume that $(f1)-(f4)$ hold, then \eqref{mainproblem} has infinitely many distinct radial solutions.
\end{theorem}
\begin{remark}
Reestablishing Theorem \ref{theorem:radialsolutions} here is not just for the sake of completeness, but aims to show that the use of the Radial Lemma due to Strauss \cite{St77} (see also \cite{Be83-1}) is not essential for obtaining radial solutions. Actually, what one really needs is the fact that the embedding $H^1_{\mathcal{O}}(\mathbb{R}^N)\hookrightarrow L^p(\mathbb{R}^N)$ is compact for all $2<p<2N/(N-2)$, see Corollary \ref{corollary:compactness1}, Remark \ref{remark:compactness1} and our proof developed for Theorem \ref{theorem:radialsolutions}.
\end{remark}

To state our results on the existence and multiplicity results of nonradial solutions, some notations are needed. Assume that $N\geq4$ and $2\leq M\leq N/2$. Let us fix $\tau\in \mathcal{O}(N)$ such that $\tau(x_1,x_2,x_3)=(x_2,x_1,x_3)$ for $x_1,x_2\in\mathbb{R}^M$ and $x_3\in\mathbb{R}^{N-2M}$, where $x=(x_1,x_2,x_3)\in\mathbb{R}^N=\mathbb{R}^M\times\mathbb{R}^M\times\mathbb{R}^{N-2M}$. We define
\begin{linenomath*}
\begin{equation}\label{eq:signchangingset}
  X_\tau:=\left\{u\in H^1(\mathbb{R}^N)~|~u(\tau x)=-u(x)~\text{for all}~x\in\mathbb{R}^N\right\}.
\end{equation}
\end{linenomath*}
It is clear that $X_\tau$ does not contain nontrivial radial functions. Let $H^1_{\mathcal{O}_1}(\mathbb{R}^N)$ denote the subspace of invariant functions with respect to $\mathcal{O}_1$, where $\mathcal{O}_1:=\mathcal{O}(M)\times\mathcal{O}(M)\times \text{id}\subset \mathcal{O}(N)$ acts isometrically on $H^1(\mathbb{R}^N)$. We also consider $\mathcal{O}_2:=\mathcal{O}(M)\times\mathcal{O}(M)\times\mathcal{O}(N-2M)\subset \mathcal{O}(N)$ acting isometrically on $H^1(\mathbb{R}^N)$ with the subspace of invariant functions denoted by $H^1_{\mathcal{O}_2}(\mathbb{R}^N)$. Here, we agree that the components corresponding to $N-2M$ do not exist when $N=2M$. Obviously, $H^1_{\mathcal{O}_2}(\mathbb{R}^N)$ is a subspace of $H^1_{\mathcal{O}_1}(\mathbb{R}^N)$, and $H^1_{\mathcal{O}_2}(\mathbb{R}^N)= H^1_{\mathcal{O}_1}(\mathbb{R}^N)$ when $N=2M$.

Now our results on nonradial solutions can be stated as follows.
\begin{theorem}\label{theorem:nonradialsolution}
 Assume that $(f1)-(f4)$ hold, $N\geq5$ and $N-2M\neq0$. Then \eqref{mainproblem} has a nonradial solution $v\in H^1_{\mathcal{O}_1}\cap X_\tau$ that minimizes $J$ among all the nontrivial solutions belonging to $H^1_{\mathcal{O}_1}\cap X_\tau$. In particular, $v$ changes signs and $J(v)>2c_{\text{mp}}$.
\end{theorem}
\begin{theorem}\label{theorem:nonradialsolutions}
 Assume that $(f1)-(f4)$ hold, $N=4$ or $N\geq6$, and $N-2M\neq1$. Then \eqref{mainproblem} has a nonradial solution $v_0\in H^1_{\mathcal{O}_2}\cap X_\tau$ which minimizes $J$ among all the nontrivial solutions belonging to $H^1_{\mathcal{O}_2}\cap X_\tau$, and admits infinitely many distinct nonradial solutions $\{v_k\}^{\infty}_{k=1}\subset H^1_{\mathcal{O}_2}\cap X_\tau$. In particular, all these solutions change signs, $J(v_0)>2c_{\text{mp}}$ and  $J(v_k)\to+\infty$ as $k\to \infty$.
\end{theorem}

The first paper dealing with the existence of nonradial solutions for equations of the type of \eqref{mainproblem} is due to  Bartsch and Willem \cite{Ba93}. They work in dimension
 $N=4$ and $N \geq 6$ under subcritical growth conditions and an Ambrosetti-Rabinowitz type condition. Actually the idea of considering subspaces of $H^1(\mathbb{R}^N)$ as $H^1_{\mathcal{O}_2}(\mathbb{R}^N)$ originates from \cite{Ba93}. Note also the work \cite{Lo04} in which the problem is solved when $N=5$ by introducing the $\mathcal{O}_1$ action on $H^1(\mathbb{R}^N)$. Finally in \cite{Mu12} Musso, Pacard and Wei, obtained nonradial solutions for any dimension $N \geq 2$, see also \cite{Ao16}. However in all these works stronger assumptions than $(f1)-(f4)$ need to be imposed. For example a nondegeneracy condition in \cite{Ao16,Mu12} which allows to apply a Lyapunov-Schmidt type reduction.

Let us now give some elements concerning the proofs of Theorems \ref{theorem:groundstate} - \ref{theorem:nonradialsolutions}. First note that in contrast to the approaches in \cite{Be83-1,Be83-2,Hi10,Me17} we shall work directly with the unconstrained functional $J$, in particular we shall not rely on the technique of adding one dimension of space, and use the mountain pass and symmetric mountain pass approaches. We know from \cite[Theorem A.VI]{Be83-1} that the energy functional $J$ is of class $C^1$. When one wants to show that $C^1$-functionals have critical points and when no general abstract result is available, a convenient first step is to show the existence of Palais-Smale sequences. This is usually done by using a quantitative deformation lemma (e.g., \cite[Lemma 2.3]{Wi96}) or Ekeland's variational principle \cite{Ek74}, if the functionals have certain convenient geometric structures. When considering the functional $J$ in $H^1(\mathbb{R}^N)$ or $H^1_\mathcal{O}(\mathbb{R}^N)$, we can conclude easily to the existence of Palais-Smale sequences since it is already known, in these two cases, that $J$ has both a mountain pass geometry and a symmetric mountain pass geometry, see \cite{Hi10,Je03} or Lemma \ref{lemma:geom1} below. However when trying to obtain nonradial solutions we have to consider $J$ restricted to the subspaces $H^1_{\mathcal{O}_1}(\mathbb{R}^N)\cap X_\tau$ or $H^1_{\mathcal{O}_2}(\mathbb{R}^N)\cap X_\tau$ which do not contain $H^1_{\mathcal{O}}(\mathbb{R}^N)$. Then, combined with the fact that we just assume $(f1)-(f4)$, the geometry of $J$ is not so apparent. Fortunately, inspired by \cite[Theorem 10]{Be83-2}, and at the expense of some technicalities, we manage in Lemma \ref{lemma:geom2} to justify the geometric properties which will insure the existence of Palais-Smale sequences.

We having obtained Palais-Smale sequences, the next obstacle is to show that these sequences are bounded. This step is particularly challenging under weak conditions as $(f1)-(f4)$. To overcome this obstacle we establish a new abstract result, Theorem \ref{theorem:SMPsetting}, which is based on the \emph{monotonicity trick}, in the spirit of \cite[Theorem 1.1]{Je99}. We recall that \cite[Theorem 1.1]{Je99} has been used extensively to deal with nonlinear variational partial differential equations where the existence of a bounded Palais-Smale sequence (at the mountain pass level) is problematic. Our extension, Theorem  \ref{theorem:SMPsetting}, can be used to derive \emph{multiple bounded Palais-Smale sequences} (at the symmetric mountain pass levels). Let us point out that Theorem \ref{theorem:SMPsetting} is not just an alternative to the technique of adding one  dimension of space \cite{Hi10,Je97}, which essentially works only for autonomous problems, but a tool which can be used to find multiple bounded Palais-Smale sequences of unconstrained functionals also in non-autonomous cases. Actually the derivation of Theorem \ref{theorem:SMPsetting} is, we believe, one of the interest of our paper.

Having dealt with the issue of the boundedness of Palais-Smale sequences, we must also study their convergence. Indeed since \eqref{mainproblem} is set on $\mathbb{R}^N$ and $(f1)-(f4)$ are weak conditions, some efforts are needed to show \emph{the convergence of bounded Palais-Smale sequences} (with respect to the norm of $H^1(\mathbb{R}^N)$).

The troublesome case is when we consider $J$ in $H^1(\mathbb{R}^N)$ the embedding from which into $L^p(\mathbb{R}^N)$ is not compact for any $2<p<2N/(N-2)$. To deal with this case, a usual way is to analyze the lack of compactness of bounded Palais-Smale sequences through the derivation of a decomposition result for the sequences, in the spirit of \cite{Lions85-1,Lions85-2}.   However, since the nonlinearity $f$ we consider is not assumed to have a subcritical growth of order $|t|^{p-1}$ for large $|t|$ with $2<p<2N/(N-2)$ and  $\lim_{t\to0}f(t)/t$ does not exist, we cannot use one of the many decomposition results in the literature. Fortunately, motivated by \cite[Proposition 4.2]{Ik17} and \cite[Proposition 4.4]{Me17}, we manage to establish a decomposition result only under the conditions $(f2)$ and $(f3)$, see Theorem \ref{theorem:decompositiontheorem}. With the aid of Theorem \ref{theorem:decompositiontheorem}, we can recover compactness at the mountain pass level $c_{mp}>0$ in the following sense: let $\{u_n\}\subset H^1(\mathbb{R}^N)$ be any bounded Palais-Smale sequence of $J$ at the level $c_{mp}>0$, up to a subsequence, there exists a sequence $\{y_n\}\subset\mathbb{R}^N$ such that the translated Palais-Smale sequence $\{u_n(\cdot+y_n)\}$ is strongly convergent in $H^1(\mathbb{R}^N)$, see Lemma \ref{lemma:compactness1}.

In the same spirit, when we consider the functional $J$ in $X=H^1_{\mathcal{O}_1}(\mathbb{R}^N)\cap X_\tau$ with $2\leq M<N/2$, a variant of Theorem \ref{theorem:decompositiontheorem} can be established and then, up to translations in $\{0\}\times\{0\}\times\mathbb{R}^{N-2M}$, the compactness can be regained at the mountain pass level for the restricted functional $J|_X$, see Corollary \ref{corollary:decompositiontheorem} and Lemma \ref{lemma:compactness2}. In the easier case where $J$ is considered in $X=H^1_{\mathcal{O}}(\mathbb{R}^N)$ with $N\geq3$ or $X=H^1_{\mathcal{O}_2}(\mathbb{R}^N)\cap X_\tau$ with $N\geq4$ and $N-2M\neq1$, the compactness issue can be addressed completely and satisfactorily. Indeed, since the embedding $X\hookrightarrow L^p(\mathbb{R}^N)$ is compact for all $2<p<2N/(N-2)$, in view of the proof of Theorem \ref{theorem:decompositiontheorem}, we can show that the restricted functional $J|_X$ satisfies \emph{the bounded Palais-Smale condition}, see Corollaries \ref{corollary:compactness1} and \ref{corollary:compactness2}.

This paper is organized as follows.

In Section \ref{sect:monotonicitytrick}, we establish Theorem \ref{theorem:SMPsetting} the new abstract result which is built under the symmetric mountain pass setting. In Section \ref{sect:decoposition}, merely under the conditions $(f2)$ and $(f3)$, we establish Theorem \ref{theorem:decompositiontheorem}, our decomposition result and present several variants. In Section \ref{sect:functionals}, still in preparation of the proofs of our main results, we introduce a family of $C^1$-functionals and work out several uniform geometric properties of them. In Section \ref{sect:proofs}, we complete the proofs of Theorems \ref{theorem:groundstate}-\ref{theorem:nonradialsolutions} by the mountain pass and symmetric mountain pass approaches.
In Section \ref{sect:remarks}, we make some miscellaneous but interesting remarks. For example, using Theorems \ref{theorem:nonradialsolution} and \ref{theorem:nonradialsolutions} and a scaling argument from \cite{Lu17}, we show the existence and multiplicity of nonradial sign-changing solutions for an autonomous Kirchhoff-type equation with Berestycki-Lions nonlinearity, these results being original. Finally, in the \hyperref[sect:appendix]{Appendix}, we prove Lemma \ref{lemma:divengence} which says that the sequence of symmetric mountain pass values goes to $+ \infty$. This property is used to show the multiplicity result of nonradial solutions claimed in Theorem \ref{theorem:nonradialsolutions}.

\section{Monotonicity trick}\label{sect:monotonicitytrick}
Assume that $(X,\|\cdot\|)$ is a real Banach space with dual space $X^{-1}$, $\mathcal{I}\subset(0,\infty)$ is a nonempty compact interval. Let $\{I_\lambda\}_{\lambda\in\mathcal{I}}$ be a family of $C^1$-functionals on $X$ being of the form
  \begin{linenomath*}
    \begin{equation*}
      I_\lambda(u)=A(u)-\lambda B(u)\qquad\text{for every}~\lambda\in\mathcal{I},
    \end{equation*}
  \end{linenomath*}
where $A,B$ are both functionals of class $C^1$, $A(0)=0=B(0)$, $B$ is nonnegative on $X$, and either $A(u)\to+\infty$ or $B(u)\to+\infty$ as $\|u\|\to \infty$.

\subsection{Mountain pass setting}\label{subsect:MPsetting}
We say that $\{I_\lambda\}_{\lambda\in\mathcal{I}}$ has a \emph{uniform mountain pass geometry} if, for every $\lambda\in\mathcal{I}$, the set
  \begin{linenomath*}
     \begin{equation*}
        \Gamma_\lambda:=\left\{\gamma\in C([0,1],X)~|~\gamma(0)=0,~I_\lambda(\gamma(1))<0\right\}
     \end{equation*}
  \end{linenomath*}
is nonempty and
  \begin{linenomath*}
     \begin{equation*}
         c_{mp,\lambda}:=\underset{\gamma\in\Gamma_\lambda}{\inf}\underset{t\in[0,1]}{\max}I_\lambda(\gamma(t))>0.
     \end{equation*}
  \end{linenomath*}
The following result is an alternative version of \cite[Theorem 1.1]{Je99} which is well suited to our needs.
\begin{theorem}[{\cite[Theorem 1.1]{Je99}}]\label{theorem:MPsetting}
  If $\{I_\lambda\}_{\lambda\in\mathcal{I}}$ has a uniform mountain pass geometry, then
   \begin{itemize}
     \item[$(i)$] for almost every $\lambda\in\mathcal{I}$, $I_\lambda$ admits a bounded Palais-Smale sequence $\{u^\lambda_n\}\subset X$ at the mountain pass level $c_{mp,\lambda}$, that is,
             \begin{linenomath*}
                \begin{equation*}
                   \sup_{n\in\mathbb{N}}\|u^\lambda_n\|<\infty,\qquad I_\lambda(u^\lambda_n)\to c_{mp,\lambda}\qquad\text{and}\qquad I'_{\lambda}(u^\lambda_n)\to0~\text{in}~X^{-1};
                \end{equation*}
             \end{linenomath*}
     \item[$(ii)$] the mapping $\lambda\mapsto c_{mp,\lambda}$ is left continuous.
   \end{itemize}
\end{theorem}
\proof Item $(i)$ can be proved by modifying the proof of \cite[Theorem 1.1]{Je99} accordingly, since
  \begin{linenomath*}
     \begin{equation*}
       \Gamma_{\lambda'}\subset\Gamma_\lambda\qquad\text{for all}~\lambda'<\lambda.
     \end{equation*}
  \end{linenomath*}
For any given $\gamma\in\Gamma_\lambda$, there exists $\delta=\delta(\gamma)>0$ such that $(\lambda-\delta,\lambda)\subset\mathcal{I}$ and
  \begin{linenomath*}
     \begin{equation*}
       \gamma\in\Gamma_{\lambda'}\qquad\text{for all}~\lambda'\in(\lambda-\delta,\lambda).
     \end{equation*}
  \end{linenomath*}
Thus, arguing as the proof of \cite[Lemma 2.3]{Je99}, we obtain Item $(ii)$.~~$\square$
\begin{remark}\label{remark:relatedresults1}
The idea to make uses of the monotonicity of the dependence of some minimax values upon a real parameter is due to Struwe \cite{St90} who however use it only on specific problems. The first abstract version of this trick is formulated in \cite{Je99} in a mountain pass setting. See also \cite{Je98}, where the condition that $B$ is nonnegative on $X$ is removed, at the expense of loosing the continuity from the left of the mapping $\lambda\mapsto c_{mp,\lambda}$. Finally let us mention \cite{Am08} in which the result of \cite{Je99} is rebuilt under a abstract minimax setting.
\end{remark}

\subsection{Symmetric mountain pass setting}\label{subsect:SMPsetting}
When $A,B$ are even, we can extend \cite[Theorem 1.1]{Je99}. For this purpose, we need to introduce a new geometric condition. For every $k\in\mathbb{N}$, let $\mathbb{D}_k:=\{x\in\mathbb{R}^k~|~|x|\leq1\}$ and $\mathbb{S}^{k-1}:=\{x\in\mathbb{R}^k~|~|x|=1\}$.

A family of even functionals $\{I_\lambda\}_{\lambda\in\mathcal{I}}$ is said to have a \emph{uniform symmetric mountain pass geometry} if, for every $k\in\mathbb{N}$, there exists an odd continuous mapping $\gamma_{0k}:\mathbb{S}^{k-1}\to X\setminus\{0\}$ such that
  \begin{linenomath*}
     \begin{equation*}
        \underset{l\in\mathbb{S}^{k-1}}{\max}I_\lambda(\gamma_{0k}(l))<0\qquad\text{uniformly in}~\lambda\in\mathcal{I},
     \end{equation*}
  \end{linenomath*}
the class of mappings $\Gamma_k:=\left\{\gamma\in C(\mathbb{D}_k,X)~|~\gamma~\text{is odd and}~\gamma=\gamma_{0k}~\text{on}~\mathbb{S}^{k-1}\right\}$ is nonempty, and
  \begin{linenomath*}
     \begin{equation*}
        c_{k,\lambda}:=\underset{\gamma\in\Gamma_k}{\inf}\underset{l\in\mathbb{D}_k}{\max}I_\lambda (\gamma(l))>0.
     \end{equation*}
  \end{linenomath*}

We shall see that $I_\lambda$ has a bounded Palais-Smale sequence at the level $c_{k,\lambda}$, for every $k\in\mathbb{N}$ and almost every $\lambda\in\mathcal{I}$, if this geometric condition is satisfied. Indeed, this is guaranteed by Theorem \ref{theorem:SMPsetting} below which can be seen as a natural extension of Theorem \ref{theorem:MPsetting}.

\begin{theorem}\label{theorem:SMPsetting}
  Assume in addition that $A,B$ are even. If $\{I_\lambda\}_{\lambda\in\mathcal{I}}$ has a uniform symmetric mountain pass geometry, then
   \begin{itemize}
     \item[$(i)$] for almost every $\lambda\in\mathcal{I}$, $I_\lambda$ admits a bounded Palais-Smale sequence $\{u^\lambda_{k,n}\}\subset X$ at each level $c_{k,\lambda}$ $(k\in\mathbb{N})$, that is,
           \begin{linenomath*}
              \begin{equation*}
                \sup_{n\in\mathbb{N}}\|u^\lambda_{k,n}\|<\infty,\qquad I_\lambda(u^\lambda_{k,n})\to c_{k,\lambda}\qquad\text{and}\qquad I'_{\lambda}(u^\lambda_{k,n})\to0~\text{in}~X^{-1};
              \end{equation*}
           \end{linenomath*}
     \item[$(ii)$] for every $k\in\mathbb{N}$, the mapping $\lambda\mapsto c_{k,\lambda}$ is left continuous.
   \end{itemize}
\end{theorem}
\begin{remark}\label{remark:relatedresults2}
  Theorem \ref{theorem:SMPsetting} will provide infinitely many bounded Palais-Smale sequences if the condition that $c_{k,\lambda}\to+\infty$ as $k\to \infty$ for any fixed $\lambda\in\mathcal{I}$ is assumed further. In general, this extra condition can be verified in specific applications. A result that is similar to Theorem \ref{theorem:SMPsetting} has been established in \cite{Zo01} for even functionals but in the setting of fountain theorems.
\end{remark}

\noindent
{\bf Proof of Theorem \ref{theorem:SMPsetting}.} Since $\Gamma_k$ is independent of $\lambda$ for every $k\in\mathbb{N}$ and $I_\lambda$ is even for all $\lambda\in\mathcal{I}$, Theorem \ref{theorem:SMPsetting} can be proved by modifying the proofs of \cite[Theorem 1.1 and Lemma 2.3]{Je99} accordingly.

For every $k\in\mathbb{N}$, since the mapping $\lambda\mapsto c_{k,\lambda}$ is non-increasing, the derivative of $c_{k,\lambda}$ with respect to $\lambda$, denoted by $c'_{k,\lambda}$, exists almost everywhere. Let $\mathcal{I}_k\subset\mathcal{I}$ be the set in which $c'_{k,\lambda}$ exists and define $\mathcal{J}:=\bigcap_{k\in\mathbb{N}}\mathcal{I}_k$. Obviously, $\mathcal{J}$ is independent of $k\in\mathbb{N}$ and $\mathcal{I}\setminus\mathcal{J}$ has zero measure. We first prove Item $(ii)$ in Claim 1 below. The proof of Item $(i)$ will be completed by showing that, for any $\lambda\in\mathcal{J}$, $I_\lambda$ admits a bounded Palais-Smale sequence at each level $c_{k,\lambda}$ ($k\in\mathbb{N}$), see Claim 3. For the proof of Claim 3, a key preliminary result is established in Claim 2.

\medskip
\textbf{Claim 1.} Item $(ii)$ holds. Namely, for every $k\in\mathbb{N}$, the mapping $\lambda\mapsto c_{k,\lambda}$ is left continuous.

\smallskip
\emph{Proof of Claim 1}.~~We assume by contradiction that, for some $k_0\in\mathbb{N}$, there exist $\lambda_{0}\in\mathcal{I}$ and $\{\lambda_n\}\subset\mathcal{I}$ such that $\lambda_n<\lambda_0$ for all $n\in\mathbb{N}$, $\lambda_n\to\lambda_0$ as $n\to\infty$, but
  \begin{linenomath*}
     \begin{equation*}
        c_{k_0,\lambda_0}<\underset{n\to\infty}{\lim}c_{k_0,\lambda_n}.
     \end{equation*}
  \end{linenomath*}
Let $\delta:=\lim_{n\to\infty}c_{k_0,\lambda_n}-c_{k_0,\lambda_0}>0$. By the definition of $c_{k_0,\lambda_0}$, we can find $\gamma_{0}\in\Gamma_{k_0}$ such that
  \begin{linenomath*}
    \begin{equation*}
      \underset{l\in\mathbb{D}_{k_0}}{\max}I_{\lambda_0}(\gamma_0(l))\leq c_{k_0,\lambda_0}+\frac{1}{3}\delta.
    \end{equation*}
  \end{linenomath*}
Using the fact that $I_\lambda(u)=I_{\lambda_0}(u)+(\lambda_0-\lambda)B(u)$ for all $\lambda\in\mathcal{I}$ and $u\in X$, we have
  \begin{linenomath*}
    \begin{equation*}
      \underset{l\in\mathbb{D}_{k_0}}{\max}I_{\lambda}(\gamma_0(l))\leq c_{k_0,\lambda_0}+\frac{1}{3}\delta+(\lambda_0-\lambda)\underset{l\in\mathbb{D}_{k_0}}{\max}B(\gamma_0(l))\qquad\text{for all}~\lambda<\lambda_0.
    \end{equation*}
  \end{linenomath*}
Since $\mathbb{D}_{k_0}$ is compact and $B$ is continuous in $u\in X$, it follows that $\max_{l\in\mathbb{D}_{k_0}}B(\gamma_0(l))\leq C$ for some $C>0$. Thus, for any $n\in\mathbb{N}$ sufficiently large,
  \begin{linenomath*}
    \begin{equation*}
      \underset{l\in\mathbb{D}_{k_0}}{\max}I_{\lambda_n}(\gamma_0(l))\leq c_{k_0,\lambda_0}+\frac{2}{3}\delta<c_{k_0,\lambda_n}.
    \end{equation*}
  \end{linenomath*}
We reach a contradiction, since the definition of $c_{k_0,\lambda_n}$ gives us that
  \begin{linenomath*}
    \begin{equation*}
      \underset{l\in\mathbb{D}_{k_0}}{\max}I_{\lambda_n}(\gamma_0(l))\geq c_{k_0,\lambda_n}\qquad\text{for all}~n\in\mathbb{N}.
    \end{equation*}
  \end{linenomath*}

We now turn to the proof of Item $(i)$. Since $\mathcal{I}\setminus\mathcal{J}$ has zero measure, we only need to show that, for any $\lambda\in\mathcal{J}$, $I_\lambda$ has a bounded Palais-Smale sequence at each level $c_{k,\lambda}$ ($k\in\mathbb{N}$). For this purpose, the following technical result is helpful. Assume that $\lambda\in\mathcal{J}$ fixed and  $\{\lambda_n\}\subset\mathcal{I}$ is a strictly increasing sequence such that $\lambda_n\to\lambda$ as $n\to\infty$.

\medskip
\textbf{Claim 2.} For every $k\in\mathbb{N}$, there exist a sequence of mappings $\{\gamma_{k,n}\}\subset\Gamma_k$ and a positive constant $K=K(c'_{k,\lambda})>0$ such that the following statements hold:
 \begin{itemize}
   \item[$(S1)$] $\|\gamma_{k,n}(l)\|\leq K$ if $\gamma_{k,n}(l)$ satisfies
               \begin{linenomath*}
                 \begin{equation}\label{eq:key1}
                   I_\lambda(\gamma_{k,n}(l))\geq c_{k,\lambda}-(\lambda-\lambda_n).
                 \end{equation}
               \end{linenomath*}
   \item[$(S2)$] $\max_{l\in\mathbb{D}_k}I_\lambda(\gamma_{k,n}(l))\leq c_{k,\lambda}+(-c'_{k,\lambda}+2)(\lambda-\lambda_n)$.
 \end{itemize}

\smallskip
\emph{Proof of Claim 2.}~~For every $k\in\mathbb{N}$, since $\Gamma_k$ is independent of $\lambda$, we can find a sequence of mappings $\{\gamma_{k,n}\}\subset\Gamma_k$ such that
  \begin{linenomath*}
     \begin{equation}\label{eq:key2}
       \max_{l\in\mathbb{D}_k}I_{\lambda_n}(\gamma_{k,n}(l))\leq c_{k,\lambda_n}+(\lambda-\lambda_n).
     \end{equation}
  \end{linenomath*}
We will show that, for $n\in\mathbb{N}$ sufficiently large, $\{\gamma_{k,n}\}$ satisfies $(S1)$ and $(S2)$. When $\gamma_{k,n}(l)$ satisfies \eqref{eq:key1}, we have
  \begin{linenomath*}
     \begin{equation*}
       \frac{I_{\lambda_n}(\gamma_{k,n}(l))-I_\lambda(\gamma_{k,n}(l))}{\lambda-\lambda_n}\leq\frac{c_{k,\lambda_n}+(\lambda-\lambda_n)-c_{k,\lambda}+(\lambda-\lambda_n)}{\lambda-\lambda_n}\leq\frac{c_{k,\lambda_n}-c_{k,\lambda}}{\lambda-\lambda_n}+2.
     \end{equation*}
  \end{linenomath*}
Since $c'_{k,\lambda}$ exists, there is $n(k,\lambda)\in\mathbb{N}$ such that, for all $n\geq n(k,\lambda)$,
  \begin{linenomath*}
    \begin{equation}\label{eq:key3}
      -c'_{k,\lambda}-1\leq\frac{c_{k,\lambda_n}-c_{k,\lambda}}{\lambda-\lambda_n}\leq -c'_{k,\lambda}+1,
    \end{equation}
  \end{linenomath*}
and then
  \begin{linenomath*}
     \begin{equation*}
       \frac{I_{\lambda_n}(\gamma_{k,n}(l))-I_\lambda(\gamma_{k,n}(l))}{\lambda-\lambda_n}\leq -c'_{k,\lambda}+3.
     \end{equation*}
  \end{linenomath*}
Consequently, for all $n\geq n(k,\lambda)$,
  \begin{linenomath*}
     \begin{equation*}
       B(\gamma_{k,n}(l))=\frac{I_{\lambda_n}(\gamma_{k,n}(l))-I_\lambda(\gamma_{k,n}(l))}{\lambda-\lambda_n}\leq -c'_{k,\lambda}+3,
     \end{equation*}
  \end{linenomath*}
and then, by \eqref{eq:key2},
  \begin{linenomath*}
     \begin{equation*}
       A(\gamma_{k,n}(l))=I_{\lambda_n}(\gamma_{k,n}(l))+\lambda_nB(\gamma_{k,n}(l))\leq c_{k,\lambda_n}+(\lambda-\lambda_n)+\lambda_n(-c'_{k,\lambda}+3)\leq C.
     \end{equation*}
  \end{linenomath*}
Since either $A(u)\to+\infty$ or $B(u)\to+\infty$ as $\|u\|\to\infty$, $(S1)$ follows directly from the uniform boundedness of $A(\gamma_{k,n}(l))$ and $B(\gamma_{k,n}(l))$. The proof of $(S2)$ is also not difficult. Indeed, \eqref{eq:key3} gives that
  \begin{linenomath*}
    \begin{equation}\label{eq:key4}
      c_{k,\lambda_n}\leq c_{k,\lambda}+(-c'_{k,\lambda}+1)(\lambda-\lambda_n)\qquad\text{for all}~n\geq n(k,\lambda).
    \end{equation}
  \end{linenomath*}
Using \eqref{eq:key2}, \eqref{eq:key4} and the fact that $I_{\lambda_n}(v)\geq I_{\lambda}(v)$ for all $v\in X$, we get
  \begin{linenomath*}
    \begin{equation*}
      \max_{l\in\mathbb{D}_k}I_\lambda(\gamma_{k,n}(l))\leq \max_{l\in\mathbb{D}_k}I_{\lambda_n}(\gamma_{k,n}(l))\leq c_{k,\lambda_n}+(\lambda-\lambda_n) \leq c_{k,\lambda}+(-c'_{k,\lambda}+2)(\lambda-\lambda_n).
    \end{equation*}
  \end{linenomath*}
The proof of Claim 2 is complete.

Let $\lambda\in\mathcal{J}$ fixed. For every $k\in\mathbb{N}$ and any $\alpha>0$, we define
\begin{linenomath*}
   \begin{equation*}
     N_{k,\alpha}:=\left\{u\in X~|~\|u\|\leq K+1~\text{and}~\left|I_\lambda(u)-c_{k,\lambda}\right|\leq \alpha\right\},
   \end{equation*}
\end{linenomath*}
where $K>0$ is the positive constant given in Claim 2. By the definition of $c_{k,\lambda}$ and Claim 2, it follows directly that $N_{k,\alpha}$ is nonempty.  One should also note that $N_{k,\alpha}\subset N_{k,\beta}$ for any $0<\alpha<\beta$.

\medskip
\textbf{Claim 3.} For every $k\in\mathbb{N}$ and any given $\alpha>0$, we have
\begin{linenomath*}
   \begin{equation*}
     \inf\left\{\|I'_\lambda(u)\|_{X^{-1}}~|~u\in N_{k,\alpha}\right\}=0.
   \end{equation*}
\end{linenomath*}
This indicates the existence of a bounded Palais-Smale sequence of $I_\lambda$ at each level $c_{k,\lambda}$ ($k\in\mathbb{N}$).

\smallskip
\emph{Proof of Claim 3.}~~We assume by contradiction that, for some $k\in\mathbb{N}$, there exists $\alpha>0$ such that
\begin{linenomath*}
   \begin{equation*}
     \|I'_\lambda(u)\|_{X^{-1}}\geq \alpha\qquad\text{for any}~u\in N_{k,\alpha}.
   \end{equation*}
\end{linenomath*}
Without loss of generality, we may assume further that
\begin{linenomath*}
   \begin{equation*}
     0<\alpha<\frac{1}{2}c_{k,\lambda}.
   \end{equation*}
\end{linenomath*}
Recall that $I_\lambda$ is even. A classical deformation argument then says that there exist $\varepsilon\in(0,\alpha)$ and an odd homeomorphism $\eta:X\to X$, such that
\begin{linenomath*}
  \begin{flalign}
      (i)&~\eta(u)=u~\text{if}~|I_\lambda(u)-c_{k,\lambda}|\geq \alpha,\label{eq:key5}&\\
     (ii)&~I_\lambda(\eta(u))\leq I_\lambda(u)~\text{for any}~u\in X,\label{eq:key6}&\\
    (iii)&~I_\lambda(\eta(u))\leq c_{k,\lambda}-\varepsilon~\text{for any}~u\in X~\text{that satisfies}~\|u\|\leq K~\text{and}~I_\lambda(u)\leq c_{k,\lambda}+\varepsilon.\label{eq:key7}&
  \end{flalign}
\end{linenomath*}
Let $\{\gamma_{k,n}\}\subset\Gamma_{k}$ be the sequence of mappings obtained in Claim 2. We can find sufficiently large but fixed $m\in\mathbb{N}$ such that
\begin{linenomath*}
   \begin{equation}\label{eq:key8}
     (-c'_{k,\lambda}+2)(\lambda-\lambda_m)\leq \varepsilon.
   \end{equation}
\end{linenomath*}
We now estimate $\max_{l\in\mathbb{D}_{k}}I_\lambda(\eta(\gamma_{k,m}(l)))$.
\begin{itemize}
  \item[$\bullet$] If $u=\gamma_{k,m}(l)$ satisfies $I_\lambda(u)\leq c_{k,\lambda}-(\lambda-\lambda_m)$, we know from \eqref{eq:key6} that
                       \begin{linenomath*}
                         \begin{equation}\label{eq:key9}
                            I_\lambda(\eta(u))\leq c_{k,\lambda}-(\lambda-\lambda_m).
                         \end{equation}
                       \end{linenomath*}
  \item[$\bullet$] If $u=\gamma_{k,m}(l)$ satisfies $I_\lambda(u)> c_{k,\lambda}-(\lambda-\lambda_m)$, in view of Claim 2 and \eqref{eq:key8}, we have that
                       \begin{linenomath*}
                          \begin{equation*}
                             \|u\|\leq K\qquad\text{and}\qquad I_\lambda(u)\leq c_{k,\lambda}+\varepsilon.
                          \end{equation*}
                       \end{linenomath*}
      Then, by \eqref{eq:key7}, it follows that
                       \begin{linenomath*}
                          \begin{equation}\label{eq:key10}
                             I_\lambda(\eta(u))\leq c_{k,\lambda}-\varepsilon\leq c_{k,\lambda}-(\lambda-\lambda_m).
                          \end{equation}
                       \end{linenomath*}
\end{itemize}
Thus, combining \eqref{eq:key9} and \eqref{eq:key10}, we get
\begin{linenomath*}
   \begin{equation*}
     \underset{l\in\mathbb{D}_{k}}{\max}I_\lambda(\eta(\gamma_{k,m}(l)))\leq c_{k,\lambda}-(\lambda-\lambda_m)<c_{k,\lambda}.
   \end{equation*}
\end{linenomath*}
While, since the homeomorphism $\eta$ is odd, by \eqref{eq:key5}, it is easy to verify that $\eta(\gamma_{k,m})\in\Gamma_{k}$. Then
\begin{linenomath*}
   \begin{equation*}
     \underset{l\in\mathbb{D}_{k}}{\max}I_\lambda(\eta(\gamma_{k,m}(l)))\geq c_{k,\lambda},
   \end{equation*}
\end{linenomath*}
we reach a contradiction.

\medskip
\textbf{Conclusion.} Obviously, Item $(i)$ now follows directly from Claim 3 and the fact that $\mathcal{I}\setminus\mathcal{J}$ has zero measure. Since Item $(ii)$ is already proved in Claim 1, the proof of Theorem \ref{theorem:SMPsetting} is complete.~~$\square$

\section{Decomposition of bounded Palais-Smale sequences}\label{sect:decoposition}
In this section, motivated by \cite[Proposition 4.2]{Ik17} and \cite[Proposition 4.4]{Me17}, we establish a decomposition result of bounded Palais-Smale sequences for a subcritical autonomous $C^1$-functional. Several variants are also derived in certain special cases. These results are necessary for us to recover a sufficient compactness when we try to prove the main theorems of this paper.

For future reference let us introduce

\begin{definition}[Bounded Palais-Smale condition]
Let $X$ be a real Banach space. We say that a $C^1$-functional $I : X \to \mathbb{R}$ satisfies the bounded Palais-Smale condition  if any bounded Palais-Smale sequence for $I$ converges, up to a subsequence.
\end{definition}

\subsection{Main decomposition result}\label{subsect:maindecomposition}
We work on $H^1(\mathbb{R}^N)$ with the standard norm
  \begin{linenomath*}
    \begin{equation*}
       \|\cdot\|_{H^1(\mathbb{R}^N)}:=\left(\int_{\mathbb{R}^N}|\nabla \cdot|^2+|\cdot|^2dx\right)^{1/2},
    \end{equation*}
  \end{linenomath*}
and consider a $C^1$-functional $I:H^1(\mathbb{R}^N)\to\mathbb{R}$ of the following form
  \begin{linenomath*}
    \begin{equation*}
       I(u):=\frac{1}{2}\int_{\mathbb{R}^N}|\nabla u|^2dx-\int_{\mathbb{R}^N}G(u)dx.
    \end{equation*}
  \end{linenomath*}
Here $N\geq3$, $G(t):=\int^t_0g(s)ds$ for all $t\in\mathbb{R}$ and $g$ is a continuous (but not necessarily odd) function satisfying $(f2)$ and $(f3)$. Our main decomposition result is stated as follows.

\begin{theorem}[Main decomposition result]\label{theorem:decompositiontheorem}
 Under the above assumptions let $\{u_n\}\subset H^1(\mathbb{R}^N)$ be a bounded Palais-Smale sequence for the functional $I$ at any level $\beta\in\mathbb{R}$. Then up to a subsequence of $\{u_n\}$ there exists an integer $l\in\mathbb{N}$ and, for each $1\leq k\leq l$, there is a sequence $\left\{y^k_n\right\}\subset\mathbb{R}^N$ and an element $w^k\in H^1(\mathbb{R}^N)$ such that the following statements hold:
     \begin{itemize}
        \item[$(i)$] $y^1_n=0$ for all $n\in\mathbb{N}$, and $|y^{i}_n-y^j_n|\to\infty$ as $n\to\infty$ for $1\leq i< j\leq l$.
        \item[$(ii)$] $u_n(\cdot +y^k_n)\rightharpoonup w^k$ in $H^1(\mathbb{R}^N)$ with $I'(w^k)=0$ for all $1\leq k\leq l$, and $w^k\neq0$ if $2\leq k\leq l$.
        \item[$(iii)$] $\beta=\lim_{n\to\infty}I(u_n)=\sum^l_{k=1} I(w^k)$.
        \item[$(iv)$] Let $v^l_n:=u_n-\sum^l_{k=1}w^k(\cdot-y^k_n)$ for every $n\in\mathbb{N}$. Then $\|v^l_n\|_{H^1(\mathbb{R}^N)}\to0$ as $n\to\infty$.
     \end{itemize}
\end{theorem}
\begin{remark}\label{remark:decomposition1}
   \begin{itemize}
      \item[$(i)$] The main feature of Theorem \ref{theorem:decompositiontheorem} is that it is established under the very weak conditions $(f2)$ and $(f3)$. In particular, we do not require the existence of a limit for $g(t)/t$ as $t\to0$. To the best of our knowledge, the first decomposition result with such a feature is due to Ikoma, but it is for a non-autonomous functional involving a fractional operator, see \cite[Proposition 4.2]{Ik17}.
      \item[$(ii)$] We highlight that similar decomposition results are expected to hold for non-autonomous (or autonomous) $C^1$-functionals without (or with) nonlocal terms, e.g., the functionals considered in \cite[Theorem 5.1]{Je05}, \cite[Lemma 3.4]{Li14} and \cite[Lemma 3.6]{Zh08} but under weak conditions like $(f2)$ and $(f3)$. Of course, the conclusions may be modified according to the specific problem under study; see, e.g., \cite[Lemma 3.4]{Li14}.
   \end{itemize}
\end{remark}

To prove Theorem \ref{theorem:decompositiontheorem}, we need a Brezis-Lieb type result and a variant of
\cite[Lemma I.1]{Lions84-2} stated as follows. The Brezis-Lieb type result can be obtained by using Vitali convergence theorem as the proof of \cite[Eq. $(3.11)$]{Me17}, and the variant of Lions lemma is only a slightly modified version of \cite[Lemma 3.1]{Me17}.
\begin{lemma}[Brezis-Lieb type result]\label{lemma:brezislieblemma}
  Assume that a function $\Psi:\mathbb{R}\to\mathbb{R}$ of class $C^1$ satisfies
    \begin{linenomath*}
      \begin{equation}\label{eq:BLcondition}
        \left|\Psi'(t)\right|\leq C\left(|t|+|t|^{\frac{N+2}{N-2}}\right)\qquad\text{for all}~t\in\mathbb{R},
      \end{equation}
    \end{linenomath*}
		and let $\{u_n\}\subset H^1(\mathbb{R}^N)$  be a bounded sequence that converges almost everywhere to $u\in H^1(\mathbb{R}^N)$  and such that $\lim_{n\to\infty} \Psi(u_n)$ exists. Then
    \begin{linenomath*}
      \begin{equation*}
         \underset{n\to\infty}{\lim}\int_{\mathbb{R}^N}\Psi(u_n)dx=\int_{\mathbb{R}^N}\Psi(u)dx+\underset{n\to\infty}{\lim} \int_{\mathbb{R}^N}\Psi(u_n-u)dx.
      \end{equation*}
    \end{linenomath*}
\end{lemma}
\begin{lemma}[Variant of Lions lemma]\label{lemma:lionslemma}
Assume that a sequence $\{u_n\}\subset H^1(\mathbb{R}^N)$ is bounded and
  \begin{linenomath*}
    \begin{equation}\label{eq:vanishing}
       \underset{n\to\infty}{\lim}\underset{y\in\mathbb{R}^N}{\sup}\int_{B(y,r)}|u_n|^2dx=0\qquad \text{for some}~r>0.
    \end{equation}
  \end{linenomath*}
Then
  \begin{linenomath*}
    \begin{equation*}
      \left|\int_{\mathbb{R}^N}\Psi(u_n)dx\right|\leq  \int_{\mathbb{R}^N}\left|\Psi(u_n)\right|dx\to0\qquad\text{as}~n\to\infty,
    \end{equation*}
  \end{linenomath*}
for any continuous function $\Psi:\mathbb{R}\to\mathbb{R}$ satisfying
  \begin{linenomath*}
     \begin{equation*}\label{eq:suplinear-subcritical}
       \underset{t\to0}{\lim}\frac{\Psi(t)}{t^2}= \underset{t\to\infty}{\lim}\frac{\Psi(t)}{|t|^{\frac{2N}{N-2}}}=0.
     \end{equation*}
  \end{linenomath*}
\end{lemma}

\noindent
\textbf{Proof of Theorem \ref{theorem:decompositiontheorem}.}
For the benefit of the reader, we shall divide the proof into three steps.

\medskip
\textbf{Step 1.} Let $y^1_n=0$ for all $n\in\mathbb{N}$. By the $H^1(\mathbb{R}^N)$-boundedness of $\{u_n\}$, we have that up to a subsequence $u_n(\cdot+y^1_n)\rightharpoonup w^1$ in $H^1(\mathbb{R}^N)$ for some $w^1\in H^1(\mathbb{R}^N)$, $u(\cdot+y^1_n)\to w^1$ in $L^p_{\text{loc}}(\mathbb{R}^N)$ for all $p\in[1,2^*)$, and $u_n(\cdot+y^1_n)\to w^1$ almost everywhere in $\mathbb{R}^N$. Since $g$ satisfies $(f2)$ and $(f3)$, with the aid of \cite[Compactness Lemma 2]{St77} (or \cite[Lemma A.I]{Be83-1}),  one can conclude that 
  \begin{linenomath*}
    \begin{equation*}
      \lim_{n\to\infty}\int_{\mathbb{R}^N}\big|\left[g(u_n)-g(w^1)\right]\phi\big|dx\leq \|\phi\|_{L^\infty(\mathbb{R}^N)}\lim_{n\to\infty}\int_{\text{supp}(\phi)}\big|g(u_n)-g(w^1)\big|dx=0
    \end{equation*}
  \end{linenomath*}
for any $\phi\in C^\infty_0(\mathbb{R}^N)$. Noting that $I'(u_n)\to0$, we obtain $I'(w^1)\phi=\lim_{n\to\infty}I'(u_n)\phi=0$.
Thus $I'(w^1)=0$. Without loss of generality, we may also assume that $\lim_{n\to\infty}\int_{\mathbb{R}^N}G(u_n)dx$ exists. Set $v^1_n:=u_n-w^1(\cdot-y^1_n)=u_n-w^1$ for every $n\in\mathbb{N}$. By $(f2)$ and $(f3)$, we see that $G$ satisfies \eqref{eq:BLcondition}. In view of Lemma \ref{lemma:brezislieblemma}, we have
  \begin{linenomath*}
    \begin{equation*}\label{eq:existenceoflimit1}
      \underset{n\to\infty}{\lim}\int_{\mathbb{R}^N}G(u_n)dx=\int_{\mathbb{R}^N}G(w^1)dx+\underset{n\to\infty}{\lim}\int_{\mathbb{R}^N}G(v^1_n)dx.
    \end{equation*}
  \end{linenomath*}
Clearly, this implies that
  \begin{linenomath*}
    \begin{equation*}\label{eq:energydecomposition1}
      \beta=\underset{n\to\infty}{\lim}I(u_n)=I(w^1)+\underset{n\to\infty}{\lim}I(v^1_n).
    \end{equation*}
  \end{linenomath*}

\medskip
\textbf{Step 2.} Assume that $m\geq1$ and, for each $1\leq k\leq m$, there is a sequence $\left\{y^k_n\right\}\subset\mathbb{R}^N$ and an elements $w^k\in H^1(\mathbb{R}^N)$ such that the following statements hold:
  \begin{itemize}
    \item[$(S1)$] $y^1_n=0$ for all $n\in\mathbb{N}$, and $|y^i_n-y^j_n|\to\infty$ as $n\to\infty$ for $1\leq i< j\leq m$.
    \item[$(S2)$] $u_n(\cdot +y^k_n)\rightharpoonup w^k$ in $H^1(\mathbb{R}^N)$ with $I'(w^k)=0$ for all $1\leq k\leq m$, and $w^k\neq0$ if $2\leq k\leq m$.
    \item[$(S3)$] Let $v^m_n:=u_n-\sum^m_{k=1}w^k(\cdot-y^k_n)$ for all $n\in\mathbb{N}$. We have that $\{v^m_n\}$ is bounded in $H^1(\mathbb{R}^N)$,
                     \begin{linenomath*}
                       \begin{equation}\label{eq:existenceoflimitm}
                           \underset{n\to\infty}{\lim}\int_{\mathbb{R}^N}G(v^m_n)dx~\text{exists}
                       \end{equation}
                     \end{linenomath*}
                   and
                     \begin{linenomath*}
                        \begin{equation}\label{eq:energydecompositionm}
                           \beta=\sum^m_{k=1} I(w^k)+\underset{n\to\infty}{\lim}I(v^m_n).
                        \end{equation}
                     \end{linenomath*}
  \end{itemize}
Letting
  \begin{linenomath*}
    \begin{equation*}
       \sigma^m:=\underset{n\to\infty}{\limsup}\big(\underset{y\in\mathbb{R}^N}{\sup}\int_{B(y,1)}|v^m_n|^2dx\big),
    \end{equation*}
  \end{linenomath*}
we distinguish the two cases: \emph{non-vanishing} and \emph{vanishing}.

\smallskip
$\bullet$ \emph{Non-vanishing: that is $\sigma^m>0$. Then, up to a subsequence of $\{u_n\}$, $(S1)-(S3)$ hold for $m+1$.} Actually,
up to a subsequence, there exists a sequence $\{y^{m+1}_n\}\subset\mathbb{R}^N$ such that
  \begin{linenomath*}
    \begin{equation*}
       \underset{n\to\infty}{\lim}\int_{B(y^{m+1}_n,1)}|v^m_n|^2dx>0.
    \end{equation*}
  \end{linenomath*}
Then $|y^{m+1}_n-y^k_n|\to\infty$ for every $1\leq k\leq m$ (since $v^m_n(\cdot+y^k_n)\to0$ in $L^2_\text{loc}(\mathbb{R}^N)$) and, up to a subsequence, $v^m_n(\cdot+y^{m+1}_n)\rightharpoonup w^{m+1}$ in $H^1(\mathbb{R}^N)$ for some $w^{m+1}\in H^1(\mathbb{R}^N)\setminus\{0\}$. Clearly,
  \begin{linenomath*}
    \begin{equation*}
       u_n(\cdot+y^{m+1}_n)=v^m_n(\cdot+y^{m+1}_n)+ \sum_{k=1}^m w^k(\cdot - y_n^k + y_n^{m+1}) \rightharpoonup w^{m+1}\qquad\text{in}~H^1(\mathbb{R}^N).
    \end{equation*}
  \end{linenomath*}
As in Step 1 we obtain that $I'(w^{m+1})=0$, since $\{u_n(\cdot+y^{m+1}_n)\}$ is a bounded Palais-Smale sequence of $I$. Let $v^{m+1}_n:=v^m_n(\cdot+y^{m+1}_n)-w^{m+1}$ for every $n\in\mathbb{N}$. We know from \eqref{eq:existenceoflimitm} and Lemma \ref{lemma:brezislieblemma} that
  \begin{linenomath*}
    \begin{equation*}
       \underset{n\to\infty}{\lim}\int_{\mathbb{R}^N}G(v^m_n(\cdot+y^{m+1}_n))dx =\int_{\mathbb{R}^N}G(w^{m+1})dx+\underset{n\to\infty}{\lim}\int_{\mathbb{R}^N}G(v^{m+1}_n(\cdot+y^{m+1}_n))dx.
    \end{equation*}
  \end{linenomath*}
Then, by \eqref{eq:energydecompositionm},
  \begin{linenomath*}
     \begin{equation*}
         \begin{aligned}
            \beta&=\sum^m_{k=1} I(w^k)+\underset{n\to\infty}{\lim}I(v^m_n)=\sum^m_{k=1}I(w^k)+\underset{n\to\infty}{\lim}I(v^m_n(\cdot+y^{m+1}_n))\\
            &=\sum^m_{k=1}I(w^k)+\left[I(w^{m+1})+\underset{n\to\infty}{\lim}I(v^{m+1}_n(\cdot+y^{m+1}_n))\right]=\sum^{m+1}_{k=1}I(w^k)+\underset{n\to\infty}{\lim}I(v^{m+1}_n).
         \end{aligned}
     \end{equation*}
  \end{linenomath*}
Thus, up to a subsequence of $\{u_n\}$, $(S1)-(S3)$ hold for $m+1$.

\smallskip
$\bullet$ \emph{Vanishing: that is $\sigma^m=0$. Then Theorem \ref{theorem:decompositiontheorem} holds with $l=m$.} Actually, since we have $(S1)$, $(S2)$ and \eqref{eq:energydecompositionm}, we only need to show that $\|v^m_n\|_{H^1(\mathbb{R}^N)}\to0$ as $n\to\infty$. For this purpose, we use an argument from \cite{Ik17}. Let
  \begin{linenomath*}
    \begin{equation*}\label{eq:vanishing1}
      \nu:=-\frac{1}{2}\underset{t\to0}{\limsup}~\frac{g(t)}{t}\in(0,\infty),
    \end{equation*}
  \end{linenomath*}
and define $\varphi(t):=g(t)+\nu t$ for all $t\in\mathbb{R}$. Obviously, $\varphi$ is a continuous function satisfying $(f2)$ and $(f3)$. Since $I'(u_n)\to0$, $I'(w^k)=0$ for all $1\leq k\leq m$, and $\{v^m_n\}$ is bounded in $H^1(\mathbb{R}^N)$, we have
    \begin{linenomath*}
      \begin{flalign*}
        \min\{1&,\nu\}\|v^m_n\|^2_{H^1(\mathbb{R}^N)}\leq\int_\mathbb{\mathbb{R}^N}|\nabla v^m_n|^2+\nu|v^m_n|^2dx\\
        &=I'(u_n)v^m_n+\int_{\mathbb{R}^N}\varphi(u_n)v^m_ndx-\sum^m_{k=1}\left[I'(w^k(\cdot-y^k_n))v^m_n+\int_{\mathbb{R}^N}\varphi(w^k(x-y^k_n))v^m_ndx\right]\\
        &=o_n(1)+\int_{\mathbb{R}^N}\left[\varphi(u_n)-\sum^m_{k=1}\varphi(w^k(x-y^k_n))\right]v^m_ndx=:o_n(1)+\Lambda_n.
      \end{flalign*}
  \end{linenomath*}
We shall show that $\limsup_{n\to\infty}\Lambda_n\leq0$.

For any $n\geq1$ and $M>0$, set $\Omega_{n,M}:=\{x~|~|v^m_n(x)|\geq M\}$. By H\"{o}lder's inequality, we have
  \begin{linenomath*}
    \begin{flalign}
      \int_{\Omega_{n,M}}&\left|\varphi(u_n)-\sum^m_{k=1}\varphi(w^k(x-y^k_n))\right||v^m_n|dx\nonumber\\
        &\qquad \leq\left(\|\varphi(u_n)\|_{L^{p^*}(\Omega_{n,M})}+\sum^m_{k=1}\|\varphi(w^k(\cdot-y^k_n))\|_{L^{p^*}(\Omega_{n,M})}\right)\|v^m_n\|_{L^{2^*}(\mathbb{R}^N)},\label{eq:vanishing3}
    \end{flalign}
  \end{linenomath*}
where $p^*:=2N/(N+2)<2$. Since $\{v^m_n\}$ is bounded in $L^{2^*}(\mathbb{R}^N)$, one can obtain that
  \begin{linenomath*}
    \begin{equation*}
      C_1\geq \|v^m_n\|^{2^*}_{L^{2^*}(\mathbb{R}^N)}\geq\|v^m_n\|^{2^*}_{L^{2^*}(\Omega_{n,M})}\geq M^{2^*}\cdot\text{meas}(\Omega_{n,M}),
    \end{equation*}
  \end{linenomath*}
where $C_1>0$ is independent of $n$ and $M$. In particular,
  \begin{linenomath*}
    \begin{equation}\label{eq:vanishing4}
      \sup_{n\geq1}\text{meas}(\Omega_{n,M})\to0\qquad\text{as}~M\to\infty.
    \end{equation}
		\end{linenomath*}
Also, for any $v\in H^1(\mathbb{R}^N)$, one has
  \begin{linenomath*}
    \begin{equation*}
      \|v\|^{p^*}_{L^{p^*}(\Omega_{n,M})} \leq  \left(\int_{\Omega_{n,M}}1dx\right)^{1-p^*/2}\cdot\left(\int_{\Omega_{n,M}} |v|^{p^*\cdot\frac{2}{p^*}}dx\right)^{p^*/2}= \big[\text{meas}(\Omega_{n,M})\big]^{1-p^*/2}\cdot\|v\|^{p^*}_{L^2(\mathbb{R}^N)}.
    \end{equation*}
  \end{linenomath*}
Since $\varphi$ satisfies $(f2)$ and $(f3)$, for any $\varepsilon>0$, we can find $C_\varepsilon>0$ such that
  \begin{linenomath*}
    \begin{equation*}
      |\varphi(t)|^{p^*}\leq C_\varepsilon|t|^{p^*}+\varepsilon |t|^{2^*}\qquad\text{for all}~t\in\mathbb{R}.
    \end{equation*}
  \end{linenomath*}
Thus it follows from Holder's inequality and the boundedness of $\{u_n\}$ that
  \begin{linenomath*}
    \begin{equation*}
      \begin{aligned}
\sup_{n\geq1}&\left(\|\varphi(u_n)\|^{p^*}_{L^{p^*}(\Omega_{n,M})}+\sum^m_{k=1}\|\varphi(w^k(\cdot-y^k_n))\|^{p^*}_{L^{p^*}(\Omega_{n,M})}\right)\\
        &\qquad\leq \sup_{n\geq1}\int_{\Omega_{n,M}}\left[C_\varepsilon|u_n|^{p^*}+\varepsilon|u_n|^{2^*}+\sum^{m}_{k=1}\left(C_\varepsilon|w^k(x-y^k_n)|^{p^*}+\varepsilon|w^k(x-y^k_n)|^{2^*}\right)\right]dx\\
        &\qquad\leq C_2\left\{C_\varepsilon\cdot\sup_{n\geq1}\big[\text{meas}(\Omega_{n,M})\big]^{1-p^*/2}+\varepsilon\right\},
      \end{aligned}
    \end{equation*}
  \end{linenomath*}
where $C_2>0$ is independent of $\varepsilon$, $n$ and $M$. Clearly, by \eqref{eq:vanishing3} and \eqref{eq:vanishing4}, we get
  \begin{linenomath*}
    \begin{equation*}
      \limsup_{M\to\infty}\sup_{n\geq1}\int_{\Omega_{n,M}}\left|\varphi(u_n)-\sum^m_{k=1}\varphi(w^k(x-y^k_n))\right||v^m_n|dx\leq C_1^{1/2^*}C_2\varepsilon.
    \end{equation*}
  \end{linenomath*}
Since $\varepsilon>0$ is arbitrary, we deduce that
  \begin{linenomath*}
    \begin{equation}\label{eq:vanishing5}
      \limsup_{M\to\infty}\sup_{n\geq1}\int_{\Omega_{n,M}}\left|\varphi(u_n)-\sum^m_{k=1}\varphi(w^k(x-y^k_n))\right||v^m_n|dx=0.
    \end{equation}
  \end{linenomath*}

On the other hand, denote by $\chi_{n,M}$ the characteristic function of the set $\{x~|~|v^m_n(x)|\leq M\}$. Clearly, for each $1\leq j\leq m$ and any $R>0$, one has
  \begin{linenomath*}
    \begin{equation*}
      \begin{split}
        \int_{B_R(y^j_n)}&\chi_{n,M}\left|\varphi(u_n)-\sum^m_{k=1}\varphi(w^k(x-y^k_n))\right||v^m_n|dx\\
        &=\int_{B_R(0)}\chi_{n,M}(x+y^j_n)\left|\varphi(u_n(x+y^j_n))-\varphi(w^j)-\sum_{k\neq j}\varphi(w^k(x+y^j_n-y^k_n))\right||v^m_n(x+y^j_n)|dx\\
        &\leq M\int_{B_R(0)}\left(\left|\varphi(u_n(x+y^j_n))-\varphi(w^j)\right|+\sum_{k\neq j}\left|\varphi(w^k(x+y^j_n-y^k_n))\right|\right)dx.
      \end{split}
    \end{equation*}
  \end{linenomath*}
Since $u_n(\cdot+y^j_n)\to w^j$ in $L^p_{\text{loc}}(\mathbb{R}^N)$ for all $p\in[1,2^*)$, $|y^j_n-y^k_n|\to\infty$ for each $k\neq j$, and $\varphi$ satisfies $(f2)$ and $(f3)$, we conclude from \cite[Compactness Lemma 2]{St77} (or \cite[Lemma A.I]{Be83-1}) that
  \begin{linenomath*}
    \begin{equation}\label{eq:vanishing6}
      \lim_{n\to\infty}\int_{B_R(y^j_n)}\chi_{n,M}\left|\varphi(u_n)-\sum^m_{k=1}\varphi(w^k(x-y^k_n))\right||v^m_n|dx=0\qquad\text{for each}~j~\text{and all}~R>0.
    \end{equation}
  \end{linenomath*}
Set $V_R:=\mathbb{R}^N\setminus\cup^m_{j=1}B_R(y^j_n)$. By $(f2)$ and $(f3)$, there exists $C_3>0$ such that
  \begin{linenomath*}
    \begin{equation*}
      |\varphi(t)|\leq C_3\left(|t|+|t|^{2^*-1}\right)\qquad\text{for all}~t\in\mathbb{R}.
    \end{equation*}
  \end{linenomath*}
Thus, for each $1\leq k\leq m$, we have
  \begin{linenomath*}
    \begin{flalign}
      \int_{V_R}&\chi_{n,M}|\varphi(w^k(x-y^k_n))v^m_n|dx\leq C_3\int_{V_R}\left(|w^k(x-y^k_n)|+|w^k(x-y^k_n)|^{2^*-1}\right)|v^m_n|dx\nonumber\\
      &\qquad\leq C_3\left(\|w^k(\cdot-y^k_n)\|_{L^2(V_R)}\|v^m_n\|_{L^2(V_R)}+\|w^k(\cdot-y^k_n)\|^{2^*-1}_{L^{2^*}(V_R)}\|v^m_n\|_{L^{2^*}(V_R)}\right)\nonumber\\
      &\qquad\leq C_3\left(\|w^k\|_{L^2(B^c_R(0))}\|v^m_n\|_{L^2(\mathbb{R}^N)}+\|w^k\|^{2^*-1}_{L^{2^*}(B^c_R(0))}\|v^m_n\|_{L^{2^*}(\mathbb{R}^N)}\right)=o_R(1),\label{eq:vanishing7}
    \end{flalign}
  \end{linenomath*}
where $o_R(1)\to0^+$ uniformly in $n$ and $M$ as $R\to\infty$. Similarly, one can obtain that
  \begin{linenomath*}
    \begin{equation}\label{eq:vanishing8}
      \begin{split}
        \int_{V_R}\chi_{n,M}|\varphi(u_n)|\left(\sum^m_{k=1}|w^k(x-y^k_n)|\right)dx=o_R(1).
      \end{split}
    \end{equation}
  \end{linenomath*}
The remaining task is to estimate the term $\int_{V_R}\varphi(u_n)\chi_{n,M}u_ndx$. Since $\varphi$ satisfies $(f2)$, one can find $\tau>0$ such that $\varphi(t)t\leq 0$ for all $|t|\leq\tau$. For $p\in(2,2^*)$ fixed and any $\varepsilon>0$, by $(f3)$, there exists $C_{p,\varepsilon}>0$ such that $|\varphi(t)t|\leq \varepsilon|t|^{2^*}+C_{p,\varepsilon} |t|^p$ for all $|t|\geq\tau$. Clearly,
  \begin{linenomath*}
    \begin{equation*}
      \begin{split}
        \int_{V_R}\varphi(u_n)\chi_{n,M}u_ndx&=\int_{V_R}\varphi(\chi_{n,M}u_n)\chi_{n,M}u_ndx\\
        &\leq \int_{V_R\cap\{x|~|u_n(x)|\geq\tau\}}\varphi(\chi_{n,M}u_n)\chi_{n,M}u_ndx\\
        &\leq \varepsilon\|u_n\|^{2^*}_{L^{2^*}(\mathbb{R}^N)}+C_{p,\varepsilon}\|u_n\|^p_{L^p(V_R)}.
      \end{split}
    \end{equation*}
  \end{linenomath*}
Since $\sigma^m=0$, we know from Lemma \ref{lemma:lionslemma} that $\lim_{n\to\infty}\|v^m_n\|_{L^p(\mathbb{R}^N)}=0$, and thus
  \begin{linenomath*}
    \begin{equation*}
      \begin{split}
        \limsup_{n\to\infty}\|u_n\|_{L^p(V_R)}\leq\limsup_{n\to\infty}\left(\|v^m_n\|_{L^p(V_R)}+\sum^m_{k=1}\|w^k(\cdot-y^k_n)\|_{L^p(V_R)}\right)=o_R(1).
      \end{split}
    \end{equation*}
  \end{linenomath*}
This implies that $\limsup_{n\to\infty}\int_{V_R}\varphi(u_n)\chi_{n,M}u_ndx\leq \varepsilon C_4 +C_{p,\varepsilon} o_R(1)$, where $C_4>0$ is independent of $\varepsilon$, $n$ and $R$. Letting $R\to\infty$, we obtain
  \begin{linenomath*}
    \begin{equation}\label{eq:vanishing9}
      \limsup_{R\to\infty}\left(\limsup_{n\to\infty}\int_{V_R}\varphi(u_n)\chi_{n,M}u_ndx\right)\leq \varepsilon C_4.
    \end{equation}
  \end{linenomath*}

We can now conclude the proof. Recall that $v^m_n=u_n-\sum^m_{k=1}w^k(\cdot-y^k_n)$. By \eqref{eq:vanishing7}, \eqref{eq:vanishing8} and \eqref{eq:vanishing9}, we have
  \begin{linenomath*}
    \begin{equation*}
      \limsup_{R\to\infty}\left[\limsup_{n\to\infty}\int_{V_R}\left(\varphi(u_n)-\sum^m_{k=1}\varphi(w^k(x-y^k_n))\right)\chi_{n,M}v_n^mdx\right]\leq \varepsilon C_4.
    \end{equation*}
  \end{linenomath*}
Since $\varepsilon>0$ is arbitrary, in view also of \eqref{eq:vanishing6}, we know that
  \begin{linenomath*}
    \begin{equation}\label{eq:vanishing90}
      \limsup_{n\to\infty}\int_{\mathbb{R}^N}\left(\varphi(u_n)-\sum^m_{k=1}\varphi(w^k(x-y^k_n))\right)\chi_{n,M}v_n^mdx\leq0.
    \end{equation}
  \end{linenomath*}
Combining \eqref{eq:vanishing90} with \eqref{eq:vanishing5}, we conclude that $\limsup_{n\to\infty}\Lambda_n\leq0$. Thus $\lim_{n\to\infty}\|v^m_n\|_{H^1(\mathbb{R}^N)}=0$. The proof of the vanishing case is complete.

\medskip
\textbf{Step 3.}  We proceed by iteration as in Step 2.  Since there is a uniformly positive constant $\rho>0$ such that $\|w\|_{H^1(\mathbb{R}^N)}\geq \rho$ for any nontrivial critical point $w$ of $I$ (see, e.g., \cite[Remark 1.3]{Je03}), we have
  \begin{linenomath*}
     \begin{equation*} \underset{n\to\infty}{\lim}\|v^m_n\|^2_{H^1(\mathbb{R}^N)}=\underset{n\to\infty}{\lim}\|u_n\|^2_{H^1(\mathbb{R}^N)}-\sum^{m}_{k=1}\|w^k\|^2_{H^1(\mathbb{R}^N)}\leq \underset{n\to\infty}{\lim}\|u_n\|^2_{H^1(\mathbb{R}^N)}-(m-1)\rho^2
     \end{equation*}
  \end{linenomath*}
if $\sigma^m>0$. Thus the vanishing case must occur for some $m_0\in\mathbb{N}$ and Theorem \ref{theorem:decompositiontheorem} holds with $l=m_0$. The proof of Theorem \ref{theorem:decompositiontheorem} is complete.~~$\square$

\subsection{Some variants of Theorem \ref{theorem:decompositiontheorem}}\label{subsect:variants}
The first special case of Theorem \ref{theorem:decompositiontheorem} occurs when we consider $I$ on $H^1_{\mathcal{O}}(\mathbb{R}^N)$, the subspace of radially symmetric functions of $H^1(\mathbb{R}^N)$. In this case, we have the following compactness result.
\begin{corollary}\label{corollary:compactness1}
Each bounded Palais-Smale sequence $\{u_n\}$ of the restricted functional $I|_{H^1_{\mathcal{O}}(\mathbb{R}^N)}$ has a strongly convergent subsequence in $H^1_{\mathcal{O}}(\mathbb{R}^N)$.
\end{corollary}
\proof The proof is nothing but a direct consequence of the proof of Theorem \ref{theorem:decompositiontheorem}. Indeed, an inspection of the proof of the vanishing case in Step 2 tells us that we only need to show that $\|v^1_n\|_{L^p(\mathbb{R}^N)}\to0$ for some $p\in(2,2^*)$. This is clearly satisfied, since the embedding $H^1_{\mathcal{O}}(\mathbb{R}^N)\hookrightarrow L^p(\mathbb{R}^N)$ is compact and $v^1_n\rightharpoonup0$ in $H^1_\mathcal{O}(\mathbb{R}^N)$.~~$\square$

Now we require additionally that the function $g$ is odd. Then, the functional $I$ is well-defined on the subspaces of $X_\tau$; see \eqref{eq:signchangingset} for the definition of $X_\tau$. Under this additional assumption two variants of Theorem \ref{theorem:decompositiontheorem} which will be used to find nonradial solutions are now presented.

Assume that $2\leq M<N/2$. For a bounded sequence $\{u_n\}\subset H^1_{\mathcal{O}_1}(\mathbb{R}^N)$, we know from \cite[Corollary 3.2]{Me17} that the conclusion of Lemma \ref{lemma:lionslemma} still holds when the condition \eqref{eq:vanishing} is replaced by the following one
   \begin{linenomath*}
     \begin{equation*}
        \underset{r\to\infty}{\lim}\big(\underset{n\to\infty}{\lim}\underset{z\in\mathbb{R}^{N-2M}}{\sup}\int_{B((0,0,z),r)}|u_n|^2dx\big)=0.
     \end{equation*}
   \end{linenomath*}
Thus, redefining $\sigma^m$ (introduced in the proof of Theorem \ref{theorem:decompositiontheorem} for $m\in\mathbb{N}$) as follows
  \begin{linenomath*}
     \begin{equation*}
        \sigma^m:=\underset{r\to\infty}{\lim}\big(\underset{n\to\infty}{\lim}\underset{z\in\mathbb{R}^{N-2M}}{\sup}\int_{B((0,0,z),r)}|v^m_n|^2dx\big)
     \end{equation*}
  \end{linenomath*}
and modifying the proof of Theorem \ref{theorem:decompositiontheorem} accordingly, we  obtain the following variant of Theorem
\ref{theorem:decompositiontheorem}.
\begin{corollary}\label{corollary:decompositiontheorem}
Assume that $2\leq M<N/2$. Let $\{u_n\}\subset X_1:=H^1_{\mathcal{O}_1}(\mathbb{R}^N)\cap X_\tau$ be a bounded Palais-Smale sequence of the restricted functional $I|_{X_1}$ at any level $\beta\in\mathbb{R}$. Then up to a subsequence of $\{u_n\}$ there exists an integer $l\in\mathbb{N}$ and, for each $1\leq k\leq l$, there is a sequence $\left\{y^k_n\right\}\subset\{0\}\times\{0\}\times\mathbb{R}^{N-2M}$ and an element $w^k\in X_1$ such that the following statements hold:
  \begin{itemize}
    \item[$(i)$] $y^1_n=0$ for all $n\in\mathbb{N}$, and $|y^i_n-y^j_n|\to\infty$ as $n\to\infty$ for $1\leq i< j\leq l$.
    \item[$(ii)$] $u_n(\cdot +y^k_n)\rightharpoonup w^k$ in $X_1$ with $(I|_{X_1})'(w^k)=0$ for all $1\leq k\leq l$, and $w^k\neq0$ if $2\leq k\leq l$.
    \item[$(iii)$] $\beta=\lim_{n\to\infty}I(u_n)= \sum^l_{k=1} I(w^k)$.
    \item[$(iv)$] Let $v^l_n:=u_n-\sum^l_{k=1}w^k(\cdot-y^k_n)$ for every $n\in\mathbb{N}$. Then $\|v^l_n\|_{H^1(\mathbb{R}^N)}\to0$ as $n\to\infty$.
  \end{itemize}
\end{corollary}

When $N\geq4$ and $N-2M\neq1$, we can say more if we choose $X_2:=H^1_{\mathcal{O}_2}(\mathbb{R}^N)\cap X_\tau$ as the working space and restrict the functional $I$ to $X_2$. Indeed, in this case, the embedding $X_2\hookrightarrow L^p(\mathbb{R}^N)$ is compact for all $2<p<2N/(N-2)$, see, e.g., \cite[Corollary 1.25]{Wi96}; then, by repeating the proof of Corollary \ref{corollary:compactness1}, we can show that
\begin{corollary}\label{corollary:compactness2}
  Assume that $N\geq4$ and $N-2M\neq1$. Then each bounded Palais-Smale sequence of the restricted functional $I|_{X_2}$ has a strongly convergent subsequence in $X_2$. Namely $I|_{X_2}$ satisfies the bounded Palais-Smale condition.
\end{corollary}

\begin{remark}\label{remark:compactness1}
Note that, in Corollaries \ref{corollary:compactness1} and \ref{corollary:compactness2}, the strong convergence is directly obtained using the fact that the working space is compactly embedded into $L^p(\mathbb{R}^N)$ for all $2<p<2N/(N-2)$. It does not involve the Radial Lemma due to Strauss \cite{St77} (see also \cite{Be83-1}).
\end{remark}

\section{Approximate functionals}\label{sect:functionals}
In this section, we introduce a family of approximate functionals and work out several uniform geometric properties of them. As one will see, with suitable choices of the working spaces, the introduced family of $C^1$-functionals satisfies the assumptions of Theorems \ref{theorem:MPsetting} and \ref{theorem:SMPsetting}. This makes possible to find bounded Palais-Smale sequences of $J$.

Let
\begin{equation}\label{eq:halfoflimsup}
  \mu:=-\frac{1}{2}\underset{t\to0}{\limsup}~\frac{f(t)}{t}\in(0,\infty).
\end{equation}
We define continuous functions $f_i$, $F_i$ $(i=1,2)$ on $\mathbb{R}$ as follows:
\begin{linenomath*}
\begin{equation*}
f_1(t):=\left\{
\begin{aligned}
&\max\{f(t)+2\mu t,0\}&\text{for}~t\geq0,\\
&\min\{f(t)+2\mu t,0\}&\text{for}~t<0,
\end{aligned}
\right.
\end{equation*}
\end{linenomath*}
\begin{linenomath*}
\begin{equation*}
f_2(t):=f_1(t)-f(t)\quad\text{for}~t\in\mathbb{R},
\end{equation*}
\end{linenomath*}
and
\begin{linenomath*}
\begin{equation*}
F_i(t):=\int^t_0f_i(s)ds\quad\text{for}~t\in\mathbb{R},\quad i=1,2.
\end{equation*}
\end{linenomath*}
Condition $(f4)$ says that $F_1(\zeta)-F_2(\zeta)>0$ for some $\zeta>0$. Thus, there exists $\lambda_0\in(0,1)$ such that $\lambda_0 F_1(\zeta)-F_2(\zeta)>0$. For $t\in\mathbb{R}$ and $\lambda\in[\lambda_0,1]$, let
\begin{linenomath*}
  \begin{equation*}
    f^\lambda(t):=\lambda f_1(t)-f_2(t)\qquad\text{and}\qquad F^\lambda(t):=\int^t_0f^\lambda(s)ds.
  \end{equation*}
\end{linenomath*}
We now introduce a family of even functionals of class $C^1$ as follows
 \begin{linenomath}
   \begin{equation*}
      J_\lambda(u):=\frac{1}{2}\int_{\mathbb{R}^N}|\nabla u|^2dx-\int_{\mathbb{R}^N}F^\lambda(u)dx,
   \end{equation*}
 \end{linenomath}
where $u\in H^1(\mathbb{R}^N)$ and $\lambda\in[\lambda_0,1]$. Since $F_1$ is nonnegative, $F_2(t)\geq\mu t^2$ for all $t\in\mathbb{R}$, and
 \begin{linenomath}
   \begin{equation*}
      J_\lambda(u)=\int_{\mathbb{R}^N}\frac{1}{2}|\nabla u|^2+F_2(u)dx-\lambda\int_{\mathbb{R}^N}F_1(u)dx=:A(u)-\lambda B(u),
   \end{equation*}
 \end{linenomath}
we see that $J_\lambda$ is of the form assumed at the very beginning of Section \ref{sect:monotonicitytrick}, and
\begin{linenomath*}
\begin{equation}\label{eq:bounds}
  J(u)=J_1(u)\leq J_\lambda(u)\leq J_{\lambda_0}(u)\qquad \text{for any}~u\in H^1(\mathbb{R}^N)~\text{and}~\lambda\in[\lambda_0,1].
\end{equation}
\end{linenomath*}

Next we present some uniform geometric properties of the functionals $J_\lambda$. Since $f$ and $f^{\lambda_0}$ satisfy $(f1)-(f4)$, by Lemma 2.4 in \cite{Hi10} and \eqref{eq:bounds} above, we have the following lemma.
\begin{lemma}\label{lemma:geom1}
Assume that $N\geq3$. Then the functional $J_\lambda$ satisfies the properties stated below.
\begin{itemize}
  \item[$(i)$] There exist $r_0>0$ and $\rho_0>0$ (independent of $\lambda\in[\lambda_0,1]$) such that
               \begin{linenomath*}
                  \begin{equation*}
                   \begin{split}
                    &J_\lambda(u)\geq J(u)>0\qquad\text{for all}~0<\|u\|_{H^1(\mathbb{R}^N)}\leq r_0,\\
                    &J_\lambda(u)\geq J(u)\geq \rho_0\qquad\text{for all}~\|u\|_{H^1(\mathbb{R}^N)}= r_0.
                   \end{split}
                 \end{equation*}
               \end{linenomath*}
  \item[$(ii)$] For every $k\in\mathbb{N}$, there exists an odd continuous mapping $\gamma_{0k}:\mathbb{S}^{k-1}\to H^1_{\mathcal{O}}(\mathbb{R}^N)$ independent of $\lambda\in[\lambda_0,1]$ such that
               \begin{linenomath*}
                 \begin{equation*}
                    J(\gamma_{0k}(l))\leq J_\lambda(\gamma_{0k}(l))\leq J_{\lambda_0}(\gamma_{0k}(l))<0\qquad\text{for all}~l\in\mathbb{S}^{k-1}.
                 \end{equation*}
               \end{linenomath*}
  \end{itemize}
\end{lemma}

When $N\geq4$, we have Lemma \ref{lemma:geom2} stated below. This result can be seen as a ``nonradial" version of Item $(ii)$ of Lemma \ref{lemma:geom1} and is essential when we try to get nonradial solutions of Problem \eqref{mainproblem}. Indeed it will be used to, for example, define the family of mappings $\Gamma_k$ in \eqref{eq:keyset}, one of the key ingredients for the verification of the uniform symmetric mountain pass geometry.

\begin{lemma}\label{lemma:geom2}
Assume further that $N\geq4$. Then, for every $k\in\mathbb{N}$, there exists an odd continuous mapping $\widetilde{\gamma}_{0k}:\mathbb{S}^{k-1}\to H^1_{\mathcal{O}_2}(\mathbb{R}^N)\cap X_\tau$ independent of $\lambda\in[\lambda_0,1]$ such that
  \begin{linenomath*}
     \begin{equation*}
        J(\widetilde{\gamma}_{0k}(l))\leq J_\lambda(\widetilde{\gamma}_{0k}(l))\leq J_{\lambda_0}(\widetilde{\gamma}_{0k}(l))<0\qquad\text{for all}~l\in\mathbb{S}^{k-1}.
     \end{equation*}
  \end{linenomath*}
\end{lemma}

The proof of Lemma \ref{lemma:geom2} is inspired by \cite{Be83-2}.  For every $k\in\mathbb{N}$ and any $R>2k$, we say that $u\in N_{k,R}$ if and only if $u\in H^1(\mathbb{R})$ is even, continuous, and satisfies the following properties:
\begin{itemize}
   \item[$(i)$] $-1\leq u\leq1$ on $[0,R)$ and $u=0$ on $[R,\infty)$.
  \item[$(ii)$] $u=\pm1$ on $[0,R]$ except in $p$ subintervals $I_1,I_2,\cdots,I_p$ of $[0,R]$ with $p\leq k$.
 \item[$(iii)$] For every $1\leq j\leq p$, $I_j$ has length at most one in which $u$ is affine with $|u'(r)|=2$.
\end{itemize}
Arguing as the Steps $(b)$ and $(c)$ in \cite[Proof of Theorem 10]{Be83-2}, we obtain an odd continuous mapping
\begin{linenomath*}
  \begin{equation*}
    U_k[R;\cdot]:\mathbb{S}^{k-1}\to N_{k,R},\qquad l\mapsto U_k[R;l](r).
  \end{equation*}
\end{linenomath*}
Let $\chi$ be an even cut-off function of class $C^\infty_0$ set on $\mathbb{R}$ such that $\chi(s)=1$ if $|s|\leq1$, $\chi(s)=0$ if $|s|\geq2$, and $0\leq\chi(s)\leq1$ if $1\leq|s|\leq2$. We then define
\begin{linenomath*}
  \begin{equation*}
     \chi(R;r):=
              \left\{
                     \begin{aligned}
                        &1,&&\text{for}~r\in[0,R^2+R],\\
                        &\chi(r-R^2-R+1),&&\text{for}~r\in[R^2+R,R^2+R+1],\\
                        &0,&&\text{for}~r\in[R^2+R+1,\infty).
                     \end{aligned}
              \right.
  \end{equation*}
\end{linenomath*}
Now, an odd continuous mapping $\pi_k[R;\cdot]:\mathbb{S}^{k-1}\to H^1_{\mathcal{O}_2}(\mathbb{R}^N)\cap X_\tau$ can be introduced as follows
\begin{linenomath*}
  \begin{equation}\label{eq:mapping}
    \pi_k[R;l](x):=\zeta\cdot \psi_k[R;l](|x_1|,|x_2|)\cdot\left|U_k[R;l](|x_3|)\right|.
  \end{equation}
\end{linenomath*}
Here $\zeta>0$ is given by $(f4)$, $x=(x_1,x_2,x_3)\in\mathbb{R}^M\times\mathbb{R}^M\times\mathbb{R}^{N-2M}$,
\begin{linenomath*}
  \begin{equation*}
    \psi_k[R;l](r_1,r_2):=\left\{U_k[R;l](r_1)-U_k[R;l](r_2)\right\}\chi(R;r_1)\chi(R;r_2),
  \end{equation*}
\end{linenomath*}
and we still agree that the components corresponding to $N-2M$ do not exist when $N=2M$.

Lemma \ref{lemma:choice} below provides an estimate on the mapping $\pi_k[R;\cdot]$ and plays an essential role in the proof of Lemma  \ref{lemma:geom2}. See also Remark \ref{remark:nonemptiness}.
\begin{lemma}\label{lemma:choice}
  For every $k\in\mathbb{N}$, there exists $R(k)>2k$ independent of $l\in\mathbb{S}^{k-1}$ such that
    \begin{linenomath*}
      \begin{equation}\label{eq:choice1}
         \int_{\mathbb{R}^N}F^{\lambda_0}(\pi_k[R;l](x))dx\geq 1\qquad\text{for all}~R\geq R(k)~\text{and}~l\in\mathbb{S}^{k-1}.
      \end{equation}
    \end{linenomath*}
\end{lemma}
\proof First, we claim that there exist $R_0(k)>2k$ and $C_M>0$ (independent of $R$ and $u,v$) such that, for all $R\geq R_0(k)$ and all $u,v\in N_{k,R}$,
\begin{linenomath*}
  \begin{equation}\label{eq:choice2}
     \int^R_0\int^{R^2+R+1}_{r_1}F^{\lambda_0}\left(\zeta\cdot[u(r_1)-v(r_2)]\chi(R;r_1)\chi(R;r_2)\right)r^{M-1}_2r^{M-1}_1dr_2dr_1\geq C_MF^{\lambda_0}(\zeta)R^{3M}.
  \end{equation}
\end{linenomath*}
Indeed, since $u,v\in N_{k,R}$, we have
\begin{linenomath*}
  \begin{flalign*}
        &\quad\int^R_0\int^{R^2+R+1}_{r_1}F^{\lambda_0}\left(\zeta\cdot[u(r_1)-v(r_2)]\chi(R;r_1)\chi(R;r_2)\right)r^{M-1}_2r^{M-1}_1dr_2dr_1\\
        &=\int^R_0\int^{R^2+R}_{R}F^{\lambda_0}\left(\zeta\cdot u(r_1)\right)r^{M-1}_2r^{M-1}_1dr_2dr_1\\
        &\qquad\qquad+\int^{R}_0\left(\int^R_{r_1}+\int^{R^2+R+1}_{R^2+R}\right)F^{\lambda_0}\left(\zeta\cdot[u(r_1)-v(r_2)]\chi(R;r_2)\right)r^{M-1}_2r^{M-1}_1dr_2dr_1\\
        &\geq\int^{R-k}_0\int^{R^2+R}_RF^{\lambda_0}\left(\zeta\right)r^{M-1}_2r^{M-1}_1dr_2dr_1-k\cdot\max_{|t|\leq\zeta}|F^{\lambda_0}(t)|\int^{R}_{R-1}\int^{R^2+R}_Rr^{M-1}_2r^{M-1}_1dr_2dr_1\\
        &\qquad\qquad-\max_{|t|\leq2\zeta}|F^{\lambda_0}(t)|\int^{R}_0\left(\int^R_{r_1}+\int^{R^2+R+1}_{R^2+R}\right)r^{M-1}_2r^{M-1}_1dr_2dr_1\\
        &\geq\frac{1}{M^2}F^{\lambda_0}(\zeta)R^{3M}\left[\frac{1}{2^M}-\frac{k}{R}\frac{2^{2M}\max_{|t|\leq\zeta}|F^{\lambda_0}(t)|}{F^{\lambda_0}(\zeta)}-\left(\frac{1}{R^M}+\frac{2^{2M-1}}{R^2}\right)\frac{\max_{|t|\leq2\zeta}|F^{\lambda_0}(t)|}{F^{\lambda_0}(\zeta)}\right].
  \end{flalign*}
\end{linenomath*}
Obviously, \eqref{eq:choice2} holds with $C_M:=2^{-M-1}M^{-2}$ by choosing $R_0(k)>2k$ large enough.

We now estimate the term $\int_{\mathbb{R}^N}F^{\lambda_0}(\pi_k[R;l](x))dx$ in the case when $N-2M\geq1$. Let
\begin{linenomath*}
\begin{equation*}
  \varphi_k[R;l](r_1,r_2,r_3):=\psi_k[R;l](r_1,r_2)\cdot \left|U_k[R;l](r_3)\right|
\end{equation*}
\end{linenomath*}
and $\omega:=\omega_{N-2M-1}\omega^2_{M-1}$, where $\omega_{m-1}$ denotes the surface area of the unit sphere in $\mathbb{R}^m$. We know from \eqref{eq:choice2} that, for all $R\geq R_0(k)$ and $l\in\mathbb{S}^{k-1}$,
\begin{linenomath*}
  \begin{equation}\label{eq:choice3}
    \int^{R}_0\int^{R^2+R+1}_{r_1}F^{\lambda_0}\left(\zeta\cdot\psi_k[R;l](r_1,r_2)\right)r^{M-1}_2r^{M-1}_1dr_2dr_1\geq C_M F^{\lambda_0}(\zeta)R^{3M}>0.
  \end{equation}
\end{linenomath*}
Thus, for all $R\geq R_0(k)$ and $l\in\mathbb{S}^{k-1}$,
\begin{linenomath*}
  \begin{flalign*}
      &\qquad\int_{\mathbb{R}^N}F^{\lambda_0}\left(\pi_k[R;l](x)\right)dx\\
      &=\omega\int^{R^2+R+1}_0\int^{R^2+R+1}_0\int^R_0F^{\lambda_0}\left(\zeta\cdot\varphi_k[R;l](r_1,r_2,r_3)\right)r^{N-2M-1}_3r^{M-1}_2r^{M-1}_1dr_3dr_2dr_1\\
      &=2\omega\int^{R}_0\int^{R^2+R+1}_{r_1}\int^R_0F^{\lambda_0}\left(\zeta\cdot\varphi_k[R;l](r_1,r_2,r_3)\right)r^{N-2M-1}_3r^{M-1}_2r^{M-1}_1dr_3dr_2dr_1\\
      &\geq2\omega\int^R_0\int^{R^2+R+1}_{r_1}F^{\lambda_0}\left(\zeta\cdot\psi_k[R;l](r_1,r_2)\right)r^{M-1}_2r^{M-1}_1dr_2dr_1\int^{R-k}_0r^{N-2M-1}_3dr_3\\
      &\qquad-2\omega\cdot k\max_{|t|\leq2\zeta}|F^{\lambda_0}(t)|\int^{R}_0\int^{R^2+R+1}_{r_1}r^{M-1}_2r^{M-1}_1dr_2dr_1\int^R_{R-1}r^{N-2M-1}_3dr_3\\
      &\geq\frac{2}{N-2M}\omega F^{\lambda_0}(\zeta)R^{N+M}\left[\frac{C_M}{2^{N-2M}}-\frac{1}{R}\frac{k\max_{|t|\leq2\zeta}|F^{\lambda_0}(t)|\cdot2^{N-M}}{M^2F^{\lambda_0}(\zeta)}\right].
  \end{flalign*}
\end{linenomath*}
Obviously, we have \eqref{eq:choice1} by choosing $R(k)\geq R_0(k)$ large enough. When $N-2M=0$, with the aid of \eqref{eq:choice3}, one can conclude easily that \eqref{eq:choice1} holds for a large enough $R(k)$.~~$\square$

Using Lemma \ref{lemma:choice}, we can now prove Lemma \ref{lemma:geom2}.

\medskip
\noindent
\textbf{Proof of Lemma \ref{lemma:geom2}.}~~For every $k\in\mathbb{N}$, let $R(k)$ be the positive constant given by Lemma \ref{lemma:choice}. Thus
    \begin{linenomath*}
      \begin{equation*}
         \int_{\mathbb{R}^N}F^{\lambda_0}(\pi_k[R(k);l](x))dx\geq 1\qquad\text{for all}~l\in\mathbb{S}^{k-1}.
      \end{equation*}
    \end{linenomath*}
Since $\mathbb{S}^{k-1}$ is compact, there exists $\beta_k>0$ such that
  \begin{linenomath*}
     \begin{equation*}
       \int_{\mathbb{R}^N}\left|\nabla \pi_k[R(k);l](x)\right|^2dx\leq \beta_k\qquad\text{for all}~l\in\mathbb{S}^{k-1}.
     \end{equation*}
  \end{linenomath*}
For any $l\in\mathbb{S}^{k-1}$, setting $\widetilde{\gamma}_{0k}(l)(x):=\pi_k[R(k);l](t^{-1}x)$ with $t\geq1$ undetermined, we have
  \begin{linenomath*}
    \begin{equation*}
       \begin{split}
         J_{\lambda_0}(\widetilde{\gamma}_{0k}(l))
         &=\frac{1}{2}t^{N-2}\int_{\mathbb{R}^N}\left|\nabla\pi_k[R(k);l](x)\right|^2dx-t^N\int_{\mathbb{R}^N}F^{\lambda_0}(\pi_k[R;l](x))dx\\
         &\leq t^{N-2}\left(\frac{1}{2}\beta_k-t^2\right).
       \end{split}
    \end{equation*}
  \end{linenomath*}
In view also of \eqref{eq:bounds}, we know that $\widetilde{\gamma}_{0k}$ is the desired odd continuous mapping by choosing $t=t_k\geq1$ large sufficiently.~~$\square$

\section{Proofs of the main results}\label{sect:proofs}
In this section, we prove Theorems \ref{theorem:groundstate}-\ref{theorem:nonradialsolutions} by mountain pass and symmetry mountain pass approaches.

\subsection{Proof of Theorem \ref{theorem:groundstate}}\label{subsect:groundstate}
We prove Theorem \ref{theorem:groundstate} by developing a mountain pass argument in $H^1(\mathbb{R}^N)$. For any $\lambda\in[\lambda_0,1]$, we know from Lemma \ref{lemma:geom1} that the set
  \begin{linenomath*}
     \begin{equation*}
       \Gamma_\lambda:=\left\{\gamma\in C([0,1],H^1(\mathbb{R}^N))~|~\gamma(0)=0, J_\lambda(\gamma(1))<0\right\}
     \end{equation*}
  \end{linenomath*}
is nonempty, the mountain pass level
  \begin{linenomath*}
     \begin{equation*}
         c_{mp,\lambda}:=\underset{\gamma\in\Gamma_\lambda}{\inf}\underset{t\in[0,1]}{\max}J_\lambda(\gamma(t))
     \end{equation*}
  \end{linenomath*}
is well-defined and $c_{mp,\lambda}\geq \rho_0>0$. We also define
   \begin{linenomath*}
     \begin{equation*}
        P_\lambda(u):=\frac{N-2}{2}\int_{\mathbb{R}^N}|\nabla u|^2dx-N\int_{\mathbb{R}^N}F^\lambda(u)dx,\qquad u\in H^1(\mathbb{R}^N).
     \end{equation*}
   \end{linenomath*}
When $\lambda=1$, for simplify, we denote $\Gamma_1$ and $P_1$ by $\Gamma$ and $P$ respectively. In view of \cite[Proof of Lemma 2.1]{Je03}, we have the following lemma.
\begin{lemma}\label{lemma:optimalpath}
  Assume that $\lambda\in[\lambda_0,1]$ is fixed and $w\in H^1(\mathbb{R}^N)\setminus\{0\}$ satisfies $P_\lambda(w)=0$. Then there exists $L>1$ (sufficiently large but fixed) such that the path defined by
  \begin{linenomath*}
    \begin{equation*}
      \gamma(t):=\left\{
                        \begin{aligned}
                           &~0,& & \quad t=0,& \\
                           &w(x/(Lt)),& & \quad t\in(0,1],&
                        \end{aligned}
                 \right.
    \end{equation*}
  \end{linenomath*}
  satisfies $\gamma(0)=0$, $\gamma(1/L)=w$, $\gamma\in C([0,1],H^1(\mathbb{R}^N))$, $J_\lambda(\gamma(1))<0$ and
  \begin{linenomath*}
    \begin{equation*}
      J_\lambda(\gamma(t))<J_\lambda(w)\qquad\text{for any}~t\in[0,1]\setminus\{1/L\}.
    \end{equation*}
  \end{linenomath*}
\end{lemma}

Using Lemma \ref{lemma:optimalpath} and our decomposition result Theorem \ref{theorem:decompositiontheorem}, we can establish the compactness result stated below.
\begin{lemma}\label{lemma:compactness1}
  Assume that $\lambda\in[\lambda_0,1]$ is fixed and $\{u_n\}\subset H^1(\mathbb{R}^N)$ is a bounded Palais-Smale sequence for $J_\lambda$ at the level $c_{mp,\lambda}$. Then, up to a subsequence, there exists a sequence $\{y_n\}\subset\mathbb{R}^N$ such that the translated sequence $\{u_n(\cdot+y_n)\}$ is a convergent Palais-Smale sequence of $J_\lambda$ at the level $c_{mp,\lambda}$.
\end{lemma}
\proof Obviously, for any given sequence $\{y_n\}\subset\mathbb{R}^N$, the translated sequence $\{u_n(\cdot+y_n)\}$ is still a bounded Palais-Smale sequence of $J_\lambda$ at the level $c_{mp,\lambda}$. To prove Lemma \ref{lemma:compactness1}, we only need to find a suitable sequence $\{y_n\}\subset\mathbb{R}^N$ such that $\{u_n(\cdot+y_n)\}$ is strongly convergent in $H^1(\mathbb{R}^N)$.

Let $w\in H^1(\mathbb{R}^N)$ be any nontrivial critical point of $J_\lambda$. Poho\u{z}aev identity implies, see for example
\cite[Proposition 1]{Be83-1}, that $P_\lambda(w)=0$. By Lemma \ref{lemma:optimalpath}, a continuous path $\gamma\in C([0,1], H^1(\mathbb{R}^N))$ exists such that
  \begin{linenomath*}
     \begin{equation*}
        \gamma\in\Gamma_\lambda\qquad\text{and}\qquad\underset{t\in[0,1]}{\max}J_\lambda(\gamma(t))=J_\lambda(w).
     \end{equation*}
  \end{linenomath*}
Therefore, $J_\lambda(w)\geq c_{mp,\lambda}$. Note that the decomposition result Theorem \ref{theorem:decompositiontheorem} applies here with $I=J_\lambda$ and $\beta=c_{mp,\lambda}>0$. If $l\geq3$, or $l=2$ but $w^1\neq0$, in view of Items $(ii)$ and $(iii)$ of Theorem \ref{theorem:decompositiontheorem}, we will get a contradiction as follows:
  \begin{linenomath*}
    \begin{equation*}
      c_{mp,\lambda}\geq \sum^l_{k=1}J_\lambda(w^k)\geq2c_{mp,\lambda}>c_{mp,\lambda}.
    \end{equation*}
  \end{linenomath*}
Thus, $l=1$, or $l=2$ with $w^1=0$. Now Lemma \ref{lemma:compactness1} follows directly from Theorem
\ref{theorem:decompositiontheorem} $(i)$ and $(iv)$.~~$\square$

To prove Theorem \ref{theorem:groundstate}, we also need the following two results which will again be used in the proofs of Theorems \ref{theorem:radialsolutions}-\ref{theorem:nonradialsolutions}.
\begin{lemma}\label{lemma:boundedness}
  Assume that $\{\lambda_n\}\subset[\lambda_0,1]$ and $\{u_n\}\subset H^1(\mathbb{R}^N)$. If
    \begin{linenomath*}
      \begin{equation*}
         \sup_{n\in\mathbb{N}}J_{\lambda_n}(u_n)\leq C\qquad\text{and}\qquad\inf_{n\in\mathbb{N}}P_{\lambda_n}(u_n)\geq-C
      \end{equation*}
    \end{linenomath*}
  for some $C>0$, then $\{u_n\}$ is bounded in $H^1(\mathbb{R}^N)$.
\end{lemma}
\proof Obviously,
  \begin{linenomath*}
    \begin{equation*}
      \frac{1}{N}\int_{\mathbb{R}^N}|\nabla u_n|^2dx=J_{\lambda_n}(u_n)-\frac{1}{N}P_{\lambda_n}(u_n)\leq 2C.
    \end{equation*}
  \end{linenomath*}
Thus $\|\nabla u_n\|_{L^2(\mathbb{R}^N)}$ is bounded. Since $F_2(t)\geq \mu t^2$ for all $t\in\mathbb{R}$ and there exists $C_\mu>0$ such that
  \begin{linenomath*}
    \begin{equation*}
      0\leq F_1(t)\leq \frac{1}{2}\mu |t|^2+C_\mu|t|^{2^*}\qquad\text{for all}~t\in\mathbb{R},
    \end{equation*}
  \end{linenomath*}
we have
  \begin{linenomath*}
    \begin{equation}\label{eq:boundedness}
        \frac{1}{2}\int_{\mathbb{R}^N}|\nabla u_n|^2+\mu|u_n|^2dx\leq J_{\lambda_n}(u_n)+\int_{\mathbb{R}^N}F_1(u_n)-\frac{1}{2}\mu|u_n|^2dx\leq C+C_\mu\int_{\mathbb{R}^N}|u_n|^{2^*}dx.
    \end{equation}
  \end{linenomath*}
In view of the boundedness of $\|\nabla u_n\|_{L^2(\mathbb{R}^N)}$, we know from Sobolev imbedding theorem that $\{u_n\}$ is bounded in $L^{2^*}(\mathbb{R}^N)$. Now the claim that $\{u_n\}$ is bounded in $H^1(\mathbb{R}^N)$ follows from \eqref{eq:boundedness}.~~$\square$

\begin{lemma}\label{lemma:boundedPSsequences}
  Assume that $\{\lambda_n\}\subset[\lambda_0,1)$, $X$ is any subspace of $H^1(\mathbb{R}^N)$, and $u_n\in X$ is a critical point of the restricted functional $J_{\lambda_n}|_X$ for every $n\in\mathbb{N}$. If $\lambda_n\to1$ as $n\to\infty$, $\{u_n\}$ is bounded in $H^1(\mathbb{R}^N)$ and
    \begin{linenomath*}
      \begin{equation*}
         \lim_{n\to\infty}J_{\lambda_n}(u_n)= c
      \end{equation*}
    \end{linenomath*}
  for some $c\in\mathbb{R}$, then $\{u_n\}$ is a bounded Palais-Smale sequence of $J|_X$ at the level $c$.
\end{lemma}
\proof  Since $\{u_n\}$ is bounded in $H^1(\mathbb{R}^N)$, we know that $\{\int_{\mathbb{R}^N}F_1(u_n)dx\}$ is bounded in $\mathbb{R}$ and $\{f_1(u_n)\}$ is bounded in $X^{-1}$. Noting that
 \begin{linenomath*}
   \begin{equation*}
     \lambda_n\to1 \qquad\text{and}\qquad J_{\lambda_n}(u_n)\to c,
   \end{equation*}
 \end{linenomath*}
we have that
\begin{linenomath*}
  \begin{equation*}
    \begin{split}
       J(u_n)&=J_{\lambda_n}(u_n)+(\lambda_n-1)F_1(u_n)=J_{\lambda_n}(u_n)+o_n(1)\to c\qquad\text{in}~\mathbb{R},\\
      (J|_X)'(u_n)&=(J_{\lambda_n}|_X)'(u_n)+(\lambda_n-1)f_1(u_n)=(\lambda_n-1)f_1(u_n)\to0\qquad\text{in}~X^{-1}.
    \end{split}
  \end{equation*}
\end{linenomath*}
Thus $\{u_n\}$ is a bounded Palais-Smale sequence for $J|_X$ at the level $c$.~~$\square$

\medskip
\noindent
\textbf{Proof of Theorem \ref{theorem:groundstate}.} Let $X=H^1(\mathbb{R}^N)$. By Theorem \ref{theorem:MPsetting}, a sequence $\{\lambda_n\}\subset[\lambda_0,1)$ exists such that
  \begin{itemize}
    \item[$(i)$] $\lambda_n\to1$ as $n\to\infty$,
    \item[$(ii)$] $c_{mp,\lambda_n}\to c_{mp,1}=c_{mp}$ as $n\to\infty$,
   \item[$(iii)$] $J_{\lambda_n}$ has a bounded Palais-Smale sequence at the level $c_{mp,\lambda_n}$ for every $n\in\mathbb{N}$.
  \end{itemize}
In view of Lemma \ref{lemma:compactness1}, we get a critical point $u_n$ of $J_{\lambda_n}$ with $J_{\lambda_n}(u_n)=c_{mp,\lambda_n}$. Poho\u{z}aev identity gives that $P_{\lambda_n}(u_n)=0$ for every $n\in\mathbb{N}$. Since
 \begin{linenomath*}
   \begin{equation*}
     \sup_{n\in\mathbb{N}}J_{\lambda_n}(u_n)=\sup_{n\in\mathbb{N}}c_{mp,\lambda_n}\leq c_{mp,\lambda_0},
   \end{equation*}
 \end{linenomath*}
by Lemma \ref{lemma:boundedness}, $\{u_n\}$ is bounded in $H^1(\mathbb{R}^N)$. Note that $J_{\lambda_n}(u_n)\to c_{mp}$. We conclude from Lemma \ref{lemma:boundedPSsequences} that $\{u_n\}$ is a bounded Palais-Smale sequence of $J$ at the mountain pass level $c_{mp}$. By Lemma \ref{lemma:compactness1} again, we get a solution $u\in H^1(\mathbb{R}^N)$ of Problem \eqref{mainproblem} with $J(u)=c_{mp}$.

We now show that $u$ is indeed a ground state solution. Define
\begin{linenomath*}
  \begin{equation*}
    \begin{split}
    \mathcal{S}&:=\{w\in H^1(\mathbb{R}^N)~|~J'(w)=0,~w\neq0\},\\
    \mathcal{P}&:=\{w\in H^1(\mathbb{R}^N)~|~P(w)=0,~w\neq0\}.
    \end{split}
  \end{equation*}
\end{linenomath*}
Obviously, $u\in\mathcal{S}\subset\mathcal{P}$ and then
\begin{linenomath*}
  \begin{equation*}
    c_{mp}=J(u)\geq \inf_{w\in\mathcal{S}}J(w)\geq\inf_{w\in\mathcal{P}}J(w).
  \end{equation*}
\end{linenomath*}
Lemma \ref{lemma:optimalpath} tells us that, for any $w\in\mathcal{P}$, there exists a path $\gamma\in C([0,1], H^1(\mathbb{R}^N))$ such that
  \begin{linenomath*}
     \begin{equation*}
        \gamma\in\Gamma\qquad\text{and}\qquad\underset{t\in[0,1]}{\max}J(\gamma(t))=J(w).
     \end{equation*}
  \end{linenomath*}
Therefore, we have
\begin{linenomath*}
  \begin{equation}\label{eq:leastenergy}
    c_{mp}=J(u)=\inf_{w\in\mathcal{S}}J(w)=\inf_{w\in\mathcal{P}}J(w),
  \end{equation}
\end{linenomath*}
which implies that $u$ is ground state solution of Problem \eqref{mainproblem}.

It has been proved in \cite{By09} that any ground state solution to Problem \eqref{mainproblem} has a constant sign. Since $f$ is odd, we may assume that $u\geq0$. Then $u>0$ by the strong maximum principle.~~$\square$

\begin{remark}\label{remark:groundstate}
 Actually, the positive ground state solution $u$ that we find is radially symmetric (up to a translation) and is decreasing with respect to the radial variable, see \cite{By09}.
\end{remark}

\subsection{Proof of Theorem \ref{theorem:nonradialsolution}}\label{subsect:nonradialsolution}
The proof of Theorem \ref{theorem:nonradialsolution} is similar to that of Theorem \ref{theorem:groundstate}. Assume that $2\leq M<N/2$ and let $X_1:=H^1_{\mathcal{O}_1}(\mathbb{R}^N)\cap X_\tau$. For fixed $\lambda\in[\lambda_0,1]$, we know from Lemma \ref{lemma:geom1} $(i)$ and Lemma \ref{lemma:geom2} that the set
  \begin{linenomath*}
     \begin{equation*}
       \widetilde{\Gamma}_\lambda:=\left\{\gamma\in C([0,1],X_1)~|~\gamma(0)=0, J_\lambda(\gamma(1))<0\right\}
     \end{equation*}
  \end{linenomath*}
is nonempty, the mountain pass level
  \begin{linenomath*}
     \begin{equation*}
         \widetilde{c}_{mp,\lambda}:=\underset{\gamma\in\widetilde{\Gamma}_\lambda}{\inf}\underset{t\in[0,1]}{\max}J_\lambda(\gamma(t))
     \end{equation*}
  \end{linenomath*}
is well-defined and $\widetilde{c}_{mp,\lambda}\geq \rho_0>0$. Replacing Theorem \ref{theorem:decompositiontheorem} by Corollary \ref{corollary:decompositiontheorem} and modifying the proof of Lemma \ref{lemma:compactness1} accordingly, we have the following compactness result.
\begin{lemma}\label{lemma:compactness2}
  Assume that $\lambda\in[\lambda_0,1]$ is fixed and $\{u_n\}\subset X_1$ is any bounded Palais-Smale sequence for $J_\lambda|_{X_1}$ at the level $\widetilde{c}_{mp,\lambda}$. Then, up to a subsequence, there exists a sequence $\{y_n\}\subset\{0\}\times\{0\}\times\mathbb{R}^{N-2M}$ such that the translated sequence $\{u_n(\cdot+y_n)\}$ is a convergent Palais-Smale sequence of $J_\lambda|_{X_1}$ at the level $\widetilde{c}_{mp,\lambda}$.
\end{lemma}

In order to show that the energy of the nonradial solution is strictly larger than $2c_{mp}$, we shall make use of Lemma \ref{lemma:criticalpoint} below.
\begin{lemma}\label{lemma:criticalpoint}
    Assume that $\gamma\in\Gamma$ and $t_*\in(0,1)$. If
    \begin{linenomath*}
       \begin{equation*}
         c_{mp}=J(\gamma(t_*))>J(\gamma(t))\qquad\text{for any}~t\in[0,1]\setminus\{t_*\},
       \end{equation*}
    \end{linenomath*}
  then $J'(\gamma(t_*))=0$. Namely, $\gamma(t_*)$ is a critical point of $J$ (at the level $c_{mp}$).
\end{lemma}
\proof This result is reminiscent of \cite[Lemma 5.1]{Mo15} and can be deduced from the quantitative deformation lemma of Willem \cite[Lemma 2.3]{Wi96}. Since its proof is essentially the same as of \cite[Lemma 5.1]{Mo15}, we omit the details here.~~$\square$

\medskip
\noindent
\textbf{Proof of Theorem \ref{theorem:nonradialsolution}.} Let $X=X_1$. By Theorem \ref{theorem:MPsetting}, a sequence $\{\lambda_n\}\subset[\lambda_0,1)$ exists such that
  \begin{itemize}
    \item[$(i)$] $\lambda_n\to1$ as $n\to\infty$,
    \item[$(ii)$] $\widetilde{c}_{mp,\lambda_n}\to \widetilde{c}_{mp,1}$ as $n\to\infty$,
    \item[$(iii)$] $J_{\lambda_n}|_{X_1}$ has a bounded Palais-Smale sequence at the level $\widetilde{c}_{mp,\lambda_n}$ for every $n\in\mathbb{N}$.
  \end{itemize}
In view of Lemma \ref{lemma:compactness2}, we get a critical point $u_n$ of $J_{\lambda_n}|_{X_1}$ with $J_{\lambda_n}(u_n)=\widetilde{c}_{mp,\lambda_n}$. By the Palais principle of symmetric criticality \cite{Pa79} and Poho\u{z}aev identity, we have that  $P_{\lambda_n}(u_n)=0$ for every $n\in\mathbb{N}$. Since
 \begin{linenomath*}
   \begin{equation*}
     \sup_{n\in\mathbb{N}}J_{\lambda_n}(u_n)=\sup_{n\in\mathbb{N}}\widetilde{c}_{mp,\lambda_n}\leq \widetilde{c}_{mp,\lambda_0},
   \end{equation*}
 \end{linenomath*}
by Lemma \ref{lemma:boundedness}, $\{u_n\}$ is bounded in $H^1(\mathbb{R}^N)$. Noting that $J_{\lambda_n}(u_n)\to \widetilde{c}_{mp,1}$, we know from Lemma \ref{lemma:boundedPSsequences} that $\{u_n\}$ is a bounded Palais-Smale sequence of $J|_{X_1}$ at the mountain pass level $\widetilde{c}_{mp,1}$. Then, by Lemma \ref{lemma:compactness2} and the Palais principle of symmetric criticality \cite{Pa79}, we get a solution $v\in X_1$ of Problem \eqref{mainproblem} with $J(v)=\widetilde{c}_{mp,1}$. Obviously, $v$ is nonradial and changes signs.

We now show that $v$ minimizes the functional $J$ among all the nontrivial solutions belonging to $X_1:=H^1_{\mathcal{O}_1}(\mathbb{R}^N)\cap X_\tau$. For this purpose, we define
\begin{linenomath*}
  \begin{equation*}
    \mathcal{S}_1:=\{w\in X_1~|~J'(w)=0,~w\neq0\}.
  \end{equation*}
\end{linenomath*}
Obviously, $v\in\mathcal{S}_1$ and then
\begin{linenomath*}
  \begin{equation*}
    \widetilde{c}_{mp,1}=J(v)\geq \inf_{w\in\mathcal{S}_1}J(w).
  \end{equation*}
\end{linenomath*}
For any $w\in\mathcal{S}_1$, we know from Poho\u{z}aev identity that $P(w)=0$. Then, by Lemma \ref{lemma:optimalpath}, there exists a path $\gamma\in C([0,1], X_1)$ such that
  \begin{linenomath*}
     \begin{equation*}
        \gamma\in\widetilde{\Gamma}_1\qquad\text{and}\qquad\underset{t\in[0,1]}{\max}J(\gamma(t))=J(w).
     \end{equation*}
  \end{linenomath*}
Therefore, we have
\begin{linenomath*}
  \begin{equation*}
    \widetilde{c}_{mp,1}=J(v)=\inf_{w\in\mathcal{S}_1}J(w).
  \end{equation*}
\end{linenomath*}

The remaining task is to show that $J(v)>2c_{mp}$. Let
\begin{linenomath*}
  \begin{equation*}
    \Omega_1:=\{x\in\mathbb{R}^N~|~|x_1|>|x_2|\}\qquad\text{and}\qquad\Omega_2:=\{x\in\mathbb{R}^N~|~|x_1|<|x_2|\}.
  \end{equation*}
\end{linenomath*}
Since $v\in X_1:=H^1_{\mathcal{O}_1}(\mathbb{R}^N)\cap X_\tau$, we get $\chi_{\Omega_j}v\in H^1_0(\Omega_j)\subset H^1(\mathbb{R}^N)$, $j=1,2$. Noting that $\chi_{\Omega_1}v\neq0$ and
\begin{linenomath*}
  \begin{equation*}
    0=P(v)=P(\chi_{\Omega_1}v)+P(\chi_{\Omega_2}v)=2P(\chi_{\Omega_1}v),
  \end{equation*}
\end{linenomath*}
we have $\chi_{\Omega_1}v\in\mathcal{P}$. Then, by \eqref{eq:leastenergy},
\begin{linenomath*}
  \begin{equation*}
    J(v)=J(\chi_{\Omega_1}v)+J(\chi_{\Omega_2}v)=2J(\chi_{\Omega_1}v)\geq 2\inf_{w\in\mathcal{P}}J(w)=2 c_{mp},
  \end{equation*}
\end{linenomath*}
that is $J(v)\geq2 c_{mp}$. If $J(v)=2 c_{mp}$, then $J(\chi_{\Omega_1}v)=c_{mp}$. By Lemma \ref{lemma:optimalpath}, Lemma \ref{lemma:criticalpoint} and  \eqref{eq:leastenergy} we deduce that $\chi_{\Omega_1}v$ is a ground state solution of Problem \eqref{mainproblem}. This is however impossible since any ground state solution of Problem \eqref{mainproblem} is radially symmetric (up to a translation), see \cite{By09}. Thus, $J(v)>2 c_{mp}$ and the proof of Theorem \ref{theorem:nonradialsolution} is complete.~~$\square$

\subsection{Proof of Theorems \ref{theorem:radialsolutions} and \ref{theorem:nonradialsolutions}}\label{subsect:multiplicity}
We first prove Theorem \ref{theorem:nonradialsolutions}. Assume that $N\geq4$ and $N-2M\neq1$, and let $X_2:=H^1_{\mathcal{O}_2}(\mathbb{R}^N)\cap X_\tau$. For any $\lambda\in[\lambda_0,1]$, by Corollary \ref{corollary:compactness2}, the restricted functional $J_\lambda|_{X_2}$ satisfies the bounded Palais-Smale condition.

For any $\lambda\in[\lambda_0,1]$, we know from Lemma \ref{lemma:geom1} $(i)$ and Lemma \ref{lemma:geom2} that the set
  \begin{linenomath*}
     \begin{equation*}
       \overline{\Gamma}_\lambda:=\left\{\gamma\in C([0,1],X_2)~|~\gamma(0)=0, J_\lambda(\gamma(1))<0\right\}
     \end{equation*}
  \end{linenomath*}
is nonempty, the mountain pass level
  \begin{linenomath*}
     \begin{equation*}
         \overline{c}_{mp,\lambda}:=\underset{\gamma\in\overline{\Gamma}_\lambda}{\inf}\underset{t\in[0,1]}{\max}J_\lambda(\gamma(t))
     \end{equation*}
  \end{linenomath*}
is well-defined and $\overline{c}_{mp,\lambda}\geq \rho_0>0$. Modifying the proof of Theorem \ref{theorem:nonradialsolution} accordingly, we can find a nonradial sign-changing solution $v_0\in H^1_{\mathcal{O}_2}\cap X_\tau$ of Problem \eqref{mainproblem} such that $v_0$ minimizes the functional $J$ among all the nontrivial solutions belonging to $H^1_{\mathcal{O}_2}\cap X_\tau$ and $J(v_0)>2c_{mp}$.

To complete the proof of Theorem \ref{theorem:nonradialsolutions}, we only need to show the multiplicity result. This will be done by developing a symmetric mountain pass argument in $X_2:=H^1_{\mathcal{O}_2}(\mathbb{R}^N)\cap X_\tau$. For every $k\in \mathbb{N}$, we define a family of mappings $\Gamma_k$ by
\begin{linenomath*}
\begin{equation}\label{eq:keyset}
  \Gamma_k:=\left\{\gamma\in C(\mathbb{D}_k,X_2))~|~\gamma~\text{is odd and}~\gamma=\widetilde{\gamma}_{0k}~\text{on}~\sigma\in \mathbb{S}^{k-1}\right\},
\end{equation}
\end{linenomath*}
where $\widetilde{\gamma}_{0k}$ is introduced in Lemma \ref{lemma:geom2}.
Clearly, $\Gamma_k$ is nonempty since it contains the mapping
\begin{linenomath*}
\begin{equation*}
\gamma_k(\sigma):=\left\{
\begin{aligned}
&|\sigma|\widetilde{\gamma}_{0k}\left(\frac{\sigma}{|\sigma|}\right),&~~~~&\text{for}~\sigma\in \mathbb{D}_k\setminus\{0\},&\\
&0,&~~~~&\text{for}~\sigma=0.&\\
\end{aligned}
\right.
\end{equation*}
\end{linenomath*}
We see from Lemma \ref{lemma:geom1} $(i)$ that
\begin{linenomath*}
\begin{equation*}
  \gamma(\mathbb{D}_k)\cap \left\{u\in X_2~|~\|u\|_{X_2}=r_0\right\}\neq \emptyset\qquad\text{for all}~ \gamma\in\Gamma_k.
\end{equation*}
\end{linenomath*}
Thus the symmetric mountain pass value $c_{k,\lambda}$ of $J_\lambda|_{X_2}$ defined by
 \begin{linenomath*}
 \begin{equation*}
c_{k,\lambda}:=\underset{\gamma\in \Gamma_k}{\inf}\underset{\sigma\in \mathbb{D}_k}{\max}~J_\lambda(\gamma(\sigma))
 \end{equation*}
 \end{linenomath*}
is well defined and it satisfies $c_{k,\lambda}\geq c_{k,1}\geq\rho_0>0$.

In the proof of Theorem \ref{theorem:nonradialsolutions} we shall need the following lemma

\begin{lemma}\label{lemma:divengence}
   The sequence of symmetric mountain pass values $\{c_{k,1}\}$ converges to $+ \infty$.
\end{lemma}

Lemma \ref{lemma:divengence} will be proved by a comparison argument due to  \cite{Hi10}. We delay its proof until the \hyperref[sect:appendix]{Appendix}.

Since $J_\lambda|_{X_2}$ satisfies the bounded Palais-Smale condition, see Corollary \ref{corollary:compactness2}. To prove the multiplicity result claimed in Theorem \ref{theorem:nonradialsolutions}, we only need to show that for every $k\in\mathbb{N}$ there is a bounded Palais-Smale sequence of $J|_{X_2}$ at the symmetric mountain pass level $c_{k,1}$. Theorem \ref{theorem:SMPsetting} is suitable for this purpose.

\medskip
\noindent
\textbf{End of the proof of Theorem \ref{theorem:nonradialsolutions}.} Let $X=X_2$. By Theorem \ref{theorem:SMPsetting}, a sequence $\{\lambda_n\}\subset[\lambda_0,1)$ exists such that
  \begin{itemize}
    \item[$(i)$] $\lambda_n\to1$ as $n\to\infty$,
   \item[$(ii)$] $c_{k,\lambda_n}\to c_{k,1}$ as $n\to\infty$ for every $k\in\mathbb{N}$,
  \item[$(iii)$] $J_{\lambda_n}|_{X_2}$ has a bounded Palais-Smale sequence at the level $c_{k,\lambda_n}$ for any $n,k\in\mathbb{N}$.
  \end{itemize}
Thus, for every $n\in\mathbb{N}$ and $k\in\mathbb{N}$, the restricted functional $J_{\lambda_n}|_{X_2}$ has a critical point $u_{k,n}$ such that $J_{\lambda_n}(u_{k,n})=c_{k,\lambda_n}$. Let $k\in\mathbb{N}$ fixed. The Palais principle of symmetric criticality \cite{Pa79} and Poho\u{z}aev identity give us that $P_{\lambda_n}(u_{k,n})=0$ for every $n\in\mathbb{N}$. Note that
 \begin{linenomath*}
   \begin{equation*}
     \sup_{n\in\mathbb{N}}J_{\lambda_n}(u_{k,n})=\sup_{n\in\mathbb{N}}c_{k,\lambda_n}\leq c_{k,\lambda_0}.
   \end{equation*}
 \end{linenomath*}
Then, the $H^1(\mathbb{R}^N)$-boundedness of $\{u_{k,n}\}$ follows from Lemma \ref{lemma:boundedness}. Since
\begin{linenomath*}
  \begin{equation*}
    J_{\lambda_n}(u_{k,n})=c_{k,\lambda_n}\to c_{k,1}\qquad\text{as}~n\to\infty,
  \end{equation*}
\end{linenomath*}
we conclude from Lemma \ref{lemma:boundedPSsequences} that $\{u_{k,n}\}$ is a bounded Palais-Smale sequence of $J|_{X_2}$ at the level $c_{k,1}$. This implies that the restricted functional $J|_{X_2}$ has a critical point $v_k\in X_2$ at each level $c_{k,1}$ ($k\in\mathbb{N}$). In view of Lemma \ref{lemma:divengence}, we have
\begin{linenomath*}
  \begin{equation*}
    J(v_k)=c_{k,1}\to+\infty\qquad\text{as}~k\to\infty.
  \end{equation*}
\end{linenomath*}
By the Palais principle of symmetric criticality \cite{Pa79}, we know that $\{v_k\}\subset X_2$ is actually a sequence of nontrivial solutions to Problem \eqref{mainproblem}. The proof of Theorem \ref{theorem:nonradialsolutions} is complete.~~$\square$

We end this section by sketching briefly the proof of Theorem \ref{theorem:radialsolutions}.

\medskip
\noindent
\textbf{Proof of Theorem \ref{theorem:radialsolutions}.}  For any $\lambda\in[\lambda_0,1]$, by Corollary \ref{corollary:compactness1}, the restricted functional $J_\lambda|_{H^1_{\mathcal{O}}(\mathbb{R}^N)}$ satisfies the bounded Palais-Smale condition. We see from Lemma \ref{lemma:geom1} that, for every $k\in \mathbb{N}$, the family of mappings
\begin{linenomath*}
\begin{equation*}
  \widehat{\Gamma}_k:=\left\{\gamma\in C(\mathbb{D}_k,H^1_{\mathcal{O}}(\mathbb{R}^N))~|~\gamma~\text{is odd and}~\gamma=\gamma_{0k}~\text{on}~\sigma\in \mathbb{S}^{k-1}\right\}
\end{equation*}
\end{linenomath*}
is nonempty, the symmetric mountain pass value
 \begin{linenomath*}
 \begin{equation*}
\widehat{c}_{k,\lambda}:=\underset{\gamma\in \widehat{\Gamma}_k}{\inf}\underset{\sigma\in \mathbb{D}_k}{\max}~J_\lambda(\gamma(\sigma))
 \end{equation*}
 \end{linenomath*}
is well-defined and $\widehat{c}_{k,\lambda}\geq \widehat{c}_{k,1}\geq\rho_0>0$.  In view of \cite[Sections 2 and 3]{Hi10}, we have that $\widehat{c}_{k,1}\to+\infty$ as $k\to\infty$. With the aid of Theorem \ref{theorem:SMPsetting}, repeating the argument above, one can obtain easily infinitely many radial solutions the energies of which converge to $+ \infty$.~~$\square$

\section{Some remarks}\label{sect:remarks}
\begin{remark}\label{extensionf3}
As it is known since \cite{Be83-1}, see also \cite{Me17}, it is possible to replace $(f3)$ by the condition
\begin{itemize}
  \item[$(f3)'$]
$- \infty \leq \limsup_{t \to + \infty}f(t)/ t^{\frac{N+2}{N-2}} \leq 0.$
\end{itemize}
Indeed assume that $f$ satisfies $(f1)$, $(f2)$, $(f3)'$ and $(f4)$. If $f(t)\geq0$ for all $t\geq\zeta$, then $f$ satisfies $(f1)-(f4)$ and our theorems in Section \ref{sect:introduction} can be applied directly. Otherwise, we set
\begin{linenomath*}
\begin{equation*}
\zeta_1:=\inf\{t\geq\zeta~|~f(t)=0\}\qquad\text{and}\qquad
\widetilde{f}(t):=\left\{
\begin{aligned}
&f(t),~&\text{for}~|t|\leq \zeta_1,\\
&0,&\text{for}~|t|>\zeta_1.\\
\end{aligned}
\right.
\end{equation*}
\end{linenomath*}
The strong maximum principle tells us that any solution $u$ of the following problem
\begin{linenomath*}
  \begin{equation}\label{eq:problem}
    -\Delta{u}= \widetilde{f}(u)\quad\text{in}~\mathbb{R}^N,\qquad u\in H^1(\mathbb{R}^N)
  \end{equation}
\end{linenomath*}
satisfies $|u(x)|\leq \zeta_1$ for all $x\in\mathbb{R}^N$. Noting that $\widetilde{f}$ satisfies $(f1)-(f4)$,  by our theorems in Section \ref{sect:introduction}, we obtain finite-energy radial and nonradial solutions of \eqref{eq:problem} which are actually also the ones of \eqref{mainproblem}.
\end{remark}

\begin{remark}
Recall the solutions $v\in H^1_{\mathcal{O}_1}(\mathbb{R}^N)\cap X_\tau$ and $v_0\in H^1_{\mathcal{O}_2}(\mathbb{R}^N)\cap X_\tau$ given by Theorems \ref{theorem:nonradialsolution} and \ref{theorem:nonradialsolutions} respectively. Obviously, when $N\geq6$, $M\geq2$ and $N-2M\geq2$, we have
\begin{linenomath*}
  \begin{equation*}
     J(v_0)\geq J(v)>2c_{mp}.
  \end{equation*}
\end{linenomath*}
It seems interesting to ask the following questions:
\begin{itemize}
  \item[$(i)$] Does the equality $J(v_0)= J(v)$ hold or not?
  \item[$(ii)$] Does, up to a translation in $\{0\}\times\{0\}\times\mathbb{R}^{N-2M}$, the solution $v$ belong to $H^1_{\mathcal{O}_2}(\mathbb{R}^N)\cap X_\tau$ or not?
\end{itemize}
\end{remark}

\begin{remark}\label{remark:nonemptiness}
The proofs of \cite[Theorems 1.2 and 1.3]{Me17} (note: these two theorems are the equivalent of our Theorems \ref{theorem:nonradialsolution} and
\ref{theorem:nonradialsolutions}) rely on the abstract result \cite[Theorem 2.2]{Me17}. In order to apply
\cite[Theorem 2.2]{Me17} it is necessary to check that
\begin{itemize}
  \item[$(i)$] the set $\{\beta\in\mathbb{R}:\gamma(\Phi^\beta)\geq k\}$ introduced in \cite[Eq. (2.4)]{Me17} is non void,
  \item[$(ii)$] the value $\beta_k$ defined therein by $\beta_k:=\inf\{\beta\in\mathbb{R}:\gamma(\Phi^\beta)\geq k\}$ is a real number.
\end{itemize}
These points seem to have been overlooked in \cite{Me17}. Let us show here how our Lemma \ref{lemma:choice} can be used to fill this gap. Let $g$ be a function satisfying $(f1)$ and $(f4)$, and set $G(t):=\int^t_0g(s)ds$. For each $k\geq1$, by Lemma \ref{lemma:choice}, there exists $R(k)>0$ such that $\int_{\mathbb{R}^N}G(\pi_k[R(k);l](x))dx\geq 1$ for all $l\in\mathbb{S}^{k-1}$, where $\pi_k[R;\cdot]$ is the odd continuous mapping given by \eqref{eq:mapping}. We now define an odd continuous mapping $\pi_k[\cdot]:\mathbb{S}^{k-1}\to H^1_{\mathcal{O}_2}(\mathbb{R}^N)\cap X_\tau$ as follows:
 \begin{linenomath*}
   \begin{equation*}
      \pi_k[l](x):=\pi_k[R(k);l]\left(\|\nabla\pi_k[R(k);l]\|^{2/(N-2)}_{L^2(\mathbb{R}^N)}x\right),\qquad x\in\mathbb{R}^N~~\text{and}~~ l\in\mathbb{S}^{k-1}.
   \end{equation*}
 \end{linenomath*}
Clearly, $\int_{\mathbb{R}^N}\left|\nabla\pi_k[l](x)\right|^2dx=1$ and $\int_{\mathbb{R}^N}G(\pi_k[l](x))dx>0$ for any $l\in\mathbb{S}^{k-1}$. Following the notations in \cite{Me17}, since $\mathbb{S}^{k-1}$ is compact, one can find $\beta(k)>0$ such that $\Phi(\pi_k[l])\leq \beta(k)$ for all $l\in\mathbb{S}^{k-1}$. Thus $\pi_k[\mathbb{S}^{k-1}]\subset\Phi^{\beta(k)}$ and this implies that $\gamma(\Phi^{\beta(k)})\geq\gamma(\pi_k[\mathbb{S}^{k-1}])\geq k$.
\end{remark}

\begin{remark}\label{remark:kirchhoff}
It is shown in \cite{Lu17}, mainly by scaling arguments, how to construct starting from an arbitrary solution  to Problem \eqref{mainproblem},  a solution to the following autonomous Kirchhoff-type equation
\begin{linenomath*}
\begin{equation}\label{kirchhoff}
  -\left(a+b\int_{\mathbb{R}^N}|\nabla{u}|^2dx\right)\Delta{u}= f(u)~~\text{in}~\mathbb{R}^N,\qquad u\in H^1(\mathbb{R}^N),
\end{equation}
\end{linenomath*}
where $a\geq0,b>0$ are constants, $N\geq3$ and $f$ satisfies $(f1)-(f4)$. Since Theorems \ref{theorem:nonradialsolution} and \ref{theorem:nonradialsolutions} provide new solutions, these results immediately translated, see
\cite[Remark 5.2]{Lu17}, into new results for \eqref{kirchhoff}. Indeed, denoting the corresponding energy functional by $\mathcal{G}$, that is,
\begin{linenomath*}
\begin{equation*}
  \mathcal{G}(u):=\frac{1}{2}a\int_{\mathbb{R}^N}|\nabla u|^2dx+\frac{1}{4}b\left(\int_{\mathbb{R}^N}|\nabla{u}|^2dx\right)^2-\int_{\mathbb{R}^N}F(u)dx,
\end{equation*}
\end{linenomath*}
we can derive the following Theorems \ref{theorem:nondegenerate} and \ref{theorem:degenerate}. These results extend those in \cite{Az11,Lu17,Lu16} where only radial solutions had been obtained.
\end{remark}

\begin{theorem}\label{theorem:nondegenerate}
Assume that $a>0$ is fixed, $b>0$, $N\geq4$ and $f$ satisfies $(f1)-(f4)$. Then the following statements hold.
\begin{itemize}
  \item[$(i)$] If $N=4$, then for every $k\in \mathbb{N}$ there exists a constant $b_k>0$ such that \eqref{kirchhoff} has at least $k$ distinct nonradial sign-changing solutions for any $b\in (0,b_k)$, the energies of which are all positive. Moreover $b_k\rightarrow 0$ as $k\rightarrow\infty$.
  \item[$(ii)$] If $N=5$, then there exists a constant $b_*>0$ such that \eqref{kirchhoff} has at least one nonradial sign-changing solution with positive energy and one nonradial sign-changing solution with negative energy for any $b\in(0,b_*)$.
    \item[$(iii)$] If $N\geq6$, then for every $k\in \mathbb{N}$ there exists a constant $b_k>0$ such that \eqref{kirchhoff} has at least $k$ distinct nonradial sign-changing solutions with positive energies and $k$ distinct nonradial sign-changing solutions with negative energies for any $b\in(0,b_k)$. Moreover $b_k\rightarrow 0$ as $k\rightarrow\infty$.
\end{itemize}
\end{theorem}
\begin{theorem}\label{theorem:degenerate}
Assume that $a=0$, $b>0$, $N\geq4$ and $f$ satisfies $(f1)-(f4)$. Then the following statements hold.
\begin{itemize}
  \item[$(i)$] If $N=4$, then there exists a positive sequence $\{b_k\}$ such that, in the case where $b= b_k$, \eqref{kirchhoff} has \emph{uncountably many nonradial sign-changing solutions} $\{u_{\lambda}\}_{\lambda>0}$. Moreover, $\mathcal{G}(u_\lambda)=0$ for all $\lambda>0$, $u_{\lambda}$ satisfies
      \begin{linenomath*}
      \begin{equation*}
      \|u_{\lambda}\|_{H^1(\mathbb{R}^N)}\rightarrow \infty~~\text{as}~\lambda\rightarrow0^+\qquad\text{and}\qquad\|u_{\lambda}\|_{H^1(\mathbb{R}^N)}\rightarrow 0~~\text{as}~\lambda\rightarrow+\infty,
      \end{equation*}
      \end{linenomath*}
      and $b_k\rightarrow 0$ as $k\rightarrow\infty$.
  \item[$(ii)$] If $N=5$, then \eqref{kirchhoff} has at least one nonradial sign-changing solution for any $b>0$ the energy of which is negative.
  \item[$(iii)$] If $N\geq6$, then \eqref{kirchhoff} has infinitely many nonradial sign-changing solutions $\{u_k\}$ for any $b>0$. Moreover, $\mathcal{G}(u_k)<0$ for every $k\in \mathbb{N}$ and $\mathcal{G}(u_k)\rightarrow0$ as $k\rightarrow\infty$.
\end{itemize}
\end{theorem}

\setcounter{section}{7}
\section*{Appendix}\label{sect:appendix}
\addcontentsline{toc}{section}{Appendix}
\renewcommand{\thesubsection}{A.\arabic{subsection}}
\setcounter{theorem}{0}
\renewcommand{\thetheorem}{A.\arabic{theorem}}
\setcounter{subsection}{0}
\setcounter{equation}{0}
\renewcommand{\theequation}{A.\arabic{equation}}

In this appendix, we prove Lemma \ref{lemma:divengence} by a comparison argument. Recall that $N\geq4$, $N-2M\neq1$ and $X_2:=H^1_{\mathcal{O}_2}(\mathbb{R}^N)\cap X_\tau$. As in \cite[Subsection 2.1]{Hi10}, we consider $p_0\in (1,\frac{N+2}{N-2})$ and set
\begin{linenomath*}
\begin{equation*}
\begin{split}
&h(t):=\left\{
\begin{aligned}
&\max\{\mu t+f(t),0\},~&\text{for}~t\geq0,\\
&\min\{\mu t+f(t),0\},&\text{for}~t<0,\\
\end{aligned}
\right.
\\
&\overline{h}(t):=\left\{
\begin{aligned}
&t^{p_0}\underset{0<s\leq t}{\max}\frac{h(s)}{s^{p_0}},~~~~~~~&~\text{for}~t>0,\\
&0,&\text{for}~t=0,\\
&-|t|^{p_0}\underset{t\leq s <0}{\max}\frac{h(|s|)}{|s|^{p_0}},&~\text{for}~t<0,\\
\end{aligned}
\right.\\
&\overline{H}(t):=\int^t_0\overline{h}(s)ds,
\end{split}
\end{equation*}
\end{linenomath*}
where $f$ satisfies $(f1)-(f4)$ and $\mu>0$ is given by \eqref{eq:halfoflimsup}. Lemma 2.1 and Corollary 2.2 in \cite{Hi10} tell us that the functions $h,\overline{h}$ and $\overline{H}$ satisfy the following properties.
\begin{linenomath*}
  \begin{flalign}
      (i)&~\text{There exists}~\delta_0>0~\text{such that}~h(t)=\overline{h}(t)=0~\text{for all}~ t\in[-\delta_0,\delta_0].\nonumber&\\
     (ii)&~\text{For all}~t\in\mathbb{R},~\text{we have}~\mu t^2/2+F(t)\leq \overline{H}(t).\label{eq:comparisoncondition}&\\
    (iii)&~\text{For all}~t\in\mathbb{R},~\text{we have}~0\leq (p_0+1)\overline{H}(t)\leq\overline{h}(t)t.\label{eq:ARcondition}&\\
     (iv)&~\text{The mapping}~t\mapsto \overline{h}(t)-\mu t~\text{satisfies}~(f1)-(f4)\label{eq:BLconditions}.&
  \end{flalign}
\end{linenomath*}

We introduce a comparison functional $\Phi:X_2\to\mathbb{R}$ of class $C^1$ as follows
\begin{linenomath*}
\begin{equation*}
\Phi(u):=\frac{1}{2}\int_{\mathbb{R}^N}|\nabla u|^2dx+\frac{1}{2}\mu\int_{\mathbb{R}^N}|u|^2dx-\int_{\mathbb{R}^N}\overline{H}(u)dx,\qquad u\in X_2.
\end{equation*}
\end{linenomath*}
It is not difficult to check that $\Phi$ has a symmetric mountain pass geometry and satisfies the Palais-Smale compactness condition.

Indeed, \eqref{eq:BLconditions} implies that $0\in X_2$ is a strict local minimum point of $\Phi$, see, e.g., \cite[Lemma 1.1]{Je03}. The odd continuous mapping $\widetilde{\gamma}_{0k}$ given by Lemma \ref{lemma:geom2} is still valid here, since \eqref{eq:comparisoncondition} implies that
\begin{linenomath*}
\begin{equation}\label{eq:comparisonresult}
J(u)\geq\Phi(u)\qquad\text{for all}~u\in X_2.
\end{equation}
\end{linenomath*}
The symmetric mountain pass values of $\Phi$ can be defined as follows:
 \begin{linenomath*}
 \begin{equation*}
d_k:=\underset{\gamma\in \Gamma_k}{\inf}\underset{\sigma\in \mathbb{D}_k}{\max}~\Phi(\gamma(\sigma)),
 \end{equation*}
 \end{linenomath*}
where $\Gamma_k$ is given by \eqref{eq:keyset} and $k\in\mathbb{N}$. Thanks to the global Ambrosetti-Rabinowitz condition \eqref{eq:ARcondition}, we can show in a standard way that every Palais-Smale sequence of $\Phi$ is bounded in $X_2$, then the Palais-Smale compactness condition follows directly from \eqref{eq:BLconditions} and Corollary \ref{corollary:compactness2}.

Now, arguing as the proof of \cite[Lemma 3.2]{Hi10}, we know that $d_k$ is a critical value of $\Phi$ and
\begin{linenomath*}
\begin{equation*}
  d_k\to +\infty\qquad\text{as}~k\to\infty.
\end{equation*}
\end{linenomath*}
In view of \eqref{eq:comparisonresult}, we see that $c_{k,1}\geq d_k$ for every $k\in\mathbb{N}$ and then $c_{k,1}\to+\infty$ as $k\to\infty$. The proof of Lemma \ref{lemma:divengence} is now complete.

\section*{Acknowledgment}
\addcontentsline{toc}{section}{Acknowledgment}
The authors would like to thank Professor Norihisa Ikoma (Keio University, Japan) for valuable comments which helped them to improve a weaker version of Theorem \ref{theorem:decompositiontheorem}.


\end{document}